\documentclass[11pt]{amsart}
\usepackage{amsmath,amsfonts,amsbsy,amsgen,amscd,mathrsfs,amssymb,amsthm}
\usepackage{enumerate,mathtools}
\usepackage{cases}

\usepackage{clipboard}
\newclipboard{myclipboard}

\usepackage[T1]{fontenc} 

\usepackage{fullpage}
\usepackage{url}

\usepackage{mathrsfs}

\usepackage[colorlinks=true,linkcolor=blue]{hyperref} 
\makeatletter
\renewcommand*{\eqref}[1]{%
  \hyperref[{#1}]{\textup{\tagform@{\ref*{#1}}}}%
}
\makeatother

\numberwithin{equation}{section}

\newtheorem{theorem}{Theorem}[section]
\newtheorem*{theorem*}{Theorem}

\newtheorem{corollary}[theorem]{Corollary}
\newtheorem*{corollary*}{Corollary}
\newtheorem{lemma}[theorem]{Lemma}
\theoremstyle{definition}
\newtheorem{definition}[theorem]{Definition}

\newtheorem{remark}[theorem]{Remark}
\newtheorem{example}[theorem]{Example}

\newtheorem*{claim*}{Claim}

\newcommand{\R}{\mathbb R}

\newcommand{\Id}{\operatorname{Id}}
\newcommand{\Tr}{\operatorname{Tr}}

\newcommand{\mf}{W}
\newcommand{\mfstar}{*}
\newcommand{\prob}{\rho}
\newcommand{\sta}{\mu_a}
\newcommand{\fin}{\mu_z}
\newcommand{\gen}{\nu}
\newcommand{\speedM}{\mathcal S}
\newcommand{\speed}{S}
\newcommand{\FishM}{\mathcal I}
\newcommand{\costM}{\mathcal T}
\newcommand{\Fish}{I}
\newcommand{\phaseM}{\mathcal V}
\newcommand{\phaseMnorm}{V}
\newcommand{\HamiltonM}{\mathcal O}
\newcommand{\Hamilton}{O}
\newcommand{\Cost}{\mathcal C}
\newcommand{\ACost}{\mathcal A}

\newcommand{\cost}{C}

\newcommand{\Hent}{E}
\newcommand{\HentM}{\mathcal E}
\newcommand{\vel}{v}
\newcommand{\ttime}{\tau}
\newcommand{\noise}{\sigma}
\newcommand{\phase}{\theta}

\newcommand*\diff{\mathop{}\!\mathrm{d}}
\newcommand{\genvec}{w}
\newcommand{\dd}{n}
\newcommand{\pot}{U}
\newcommand{\f}{f}
\newcommand{\F}{F}
\newcommand{\sgn}{\epsilon}
\newcommand{\Pheat}{\mathrm{P}}
\newcommand{\confn}{\mathrm{c}}
\newcommand{\nice}{\text{nice }}
\newcommand{\hamiltonian}{H}
\newcommand{\lagrangian}{L}
\newcommand{\divnab}{\nabla{\cdot}}
\newcommand{\moment}{p}
\newcommand{\domain}{\Omega}
\newcommand{\valuefn}{\mathrm{V}}
\newcommand{\pressure}{p}

\title{Matrix displacement convexity along density flows}

\author{Yair Shenfeld}
\address{Division of Applied Mathematics, Brown University, Providence, RI, USA}
\email{Yair\_Shenfeld@Brown.edu}

\begin{document}
\maketitle

\begin{abstract}
A new notion of displacement convexity on a matrix level is developed for density flows arising from mean-field games, compressible Euler equations, entropic interpolation, and semi-classical limits of non-linear Schr\"odinger equations. Matrix displacement convexity is stronger than the classical notions of displacement convexity, and its verification (formal and rigorous)  relies on  matrix differential inequalities along the density flows. The matrical nature of these differential inequalities upgrades dimensional functional inequalities to their intrinsic dimensional counterparts, thus improving on many classical results. Applications include turnpike properties, evolution variational inequalities, and entropy growth bounds, which capture the behavior of the density flows along different directions in space.
\end{abstract}

\section{Introduction}
 
The optimal decisions of agents in large populations, the lazy gas experiment of Schr\"odinger, and the flow of slender jets can all be modeled by systems of coupled partial differential equations of the form
\begin{equation}
\label{eq:flow_intro}
\begin{cases}
\partial_t\prob_t+\divnab\left(\prob_t\nabla\phase_t\right)=0,\\
\partial_t\phase_t+\frac{1}{2}|\nabla\phase_t|^2+\frac{\noise^2}{2}\frac{\Delta \prob_t^{1/2}}{\prob_t^{1/2}}+\pot_t-\mf\mfstar \prob_t-\f(\prob_t)=0,
\end{cases}
\forall ~t\in [0,\ttime].
\end{equation}
The first equation in \eqref{eq:flow_intro} is the continuity equation of $\prob_t\ge 0$, interpreted as a probability density over a domain $\domain\subset \R^{\dd}$, driven by a  gradient vector field $\nabla\phase_t$. The second equation  in \eqref{eq:flow_intro} describes the evolution of the vector field itself  via an equation for $\phase_t$, which in turn can depend on the density $\prob_t$. The boundary conditions  of \eqref{eq:flow_intro} will vary based on the model and will usually be a specification of  $(\prob_0,\prob_{\ttime})$, or $(\prob_0,\nabla\phase_0)$, or $(\prob_0,\phase_{\ttime})$, and so on. The scope of the flows \eqref{eq:flow_intro} will be recalled in Section \ref{subsec:examples}: they include \textbf{planning problems} (optimal transport, entropic interpolation, regularization of planning problems), \textbf{mean-field games}, \textbf{barotropic fluids}, and \textbf{semi-classical limits of non-linear Schr\"odinger equations}. The majority of this work will focus on flows $(\prob_t,\phase_t)$ satisfying \eqref{eq:flow_intro} under the assumptions\footnote{Some of these assumptions can in fact be relaxed, cf. Section \ref{subsec:main_results}.} that $\noise$ is real, $\pot_t$ is convex, $\mf$ is concave, and $\f$ is non-decreasing.

\subsection{Matrix displacement convexity and intrinsic dimensional functional inequalities}
\label{subsec:martix_cvx}

The discovery by McCann \cite{McCann1997} of the notion of displacement convexity has had a significant impact on probability, analysis, and geometry. Specifically, it was shown in \cite{McCann1997} that certain functionals are convex along the  optimal transport flow. A central example of such a functional is the differential entropy
\[
\Hent(t):=\int_{\domain}\log \prob_t\diff\prob_t,
\]
which was shown by McCann \cite{McCann1997} to be convex (i.e., $t\mapsto \Hent(t)$ is convex) when $(\prob_t)$ is the optimal transport flow. It was later realized that displacement convexity also holds along other density flows. For example, L\'eonard \cite{Leonardconvex} showed that $\Hent(t)$ is convex when $(\prob_t)$ is the entropic interpolation flow, and Gomes and Seneci \cite{gomes2018displacement} showed that $\Hent(t)$ is convex when $(\prob_t)$ is a first-order mean-field game flow. In a certain sense, these results generalize McCann's result (as well as the classical  convexity of entropy along heat flows) as will become clear in  Section \ref{subsec:examples}. 

One important application of displacement convexity is its usage in the definition of Ricci curvature for metric measure spaces as developed by Lott-Villani \cite{LottVillani2009} and  Sturm  \cite{Sturm2006I,Sturm2006II}. Roughly, a metric-measure space is defined to have a nonnegative Ricci curvature if the entropy is convex along optimal transport flows over this space. This notion of Ricci curvature coincides with the classical notion when the space is a Riemannian manifold. 

There is a stronger curvature condition (going beyond nonnegative Ricci curvature) which incorporates the effect of the dimension. Restricting to the flat case, this is the $\text{CD}(0,\dd)$ curvature-dimension condition of Bakry-\'Emery \cite{BE2014}. Analogous to the relation between nonnegative Ricci curvature and displacement convexity of the entropy, Erbar-Kuwada-Sturm \cite{EKS2015} showed  that the $\text{CD}(0,\dd)$ curvature-dimension condition is equivalent (under sufficient regularity) to the  concavity of the map 
\[
t\mapsto e^{-\frac{\Hent(t)}{\dd}}
\]
along the optimal transport flow (which implies the convexity of entropy along the flow).  Due to the role of dimension in this notion of convexity it will be dubbed here \emph{dimensional displacement convexity}. The natural question of whether $e^{-\frac{\Hent(t)}{\dd}}$ is concave along the entropic interpolation flow was settled by Ripani \cite{ripani2019convexity}, which thus recovered both the result of Erbar-Kuwada-Sturm on flat space, as well as the result of Costa \cite{costa1985} who showed that $e^{-\frac{\Hent(t)}{\dd}}$ is concave along the heat flow.

\subsubsection{Matrix displacement convexity}
\label{subsubsec:matrixdspcv_intdimfuinq}
The main purpose of this work is to develop and prove a new notion of \emph{matrix} displacement convexity which is stronger than  dimensional displacement convexity (and thus stronger than classical displacement convexity) along density flows of the form \eqref{eq:flow_intro}. To keep the discussion concrete, at this point matrix displacement convexity will be defined just for the entropy (but the extension is clear). Recall that the \emph{entropy production} $\speed(t)$ associated to a density flow $(\prob_t)$ is defined as
\[
\speed(t):=\partial_t\Hent(t). 
\]
In the setting of this work (and many others), there is a natural \emph{entropy production matrix} $\speedM(t)$ which can be defined so that
\[
\partial_t\Hent(t)=\speed(t)=\Tr[\speedM(t)].
\]
Indeed, a  simple calculation (cf. Lemma \ref{lem:entropy_1der}) shows that when $(\prob_t,\phase_t)$ satisfies the continuity equation, the \emph{entropy production matrix} is given by
\[
\speedM(t)=\int_{\domain} \nabla\prob_t\otimes_S\nabla\phase_t\diff x
\]
where $\otimes_S$ is the symmetric tensor product. The \emph{entropy matrix} is defined as 
\[
\HentM(t) :=\int_0^t\speedM(s)\diff s
\]
so that
\[
\Hent(t)=\Hent(0)+\Tr[\HentM(t)].
\]
\begin{definition}
\label{def:matrix_disp_cnvx_entropy}
The entropy matrix $\HentM(t)$ is  \emph{matrix displacement convex} along a flow $(\prob_t,\phase_t)$ if, for any unit vector $\genvec\in \R^{\dd}$, the function 
\[
t\mapsto e^{-\langle \genvec,\HentM(t)\genvec\rangle}
\]
is concave. 
\end{definition}
Note that if $\HentM(t)$ is matrix displacement convex then $\Hent(t)$ is dimensional displacement convex, and hence displacement convex (cf. Section \ref{subsubsec:matrix_diff}.) 

There are two main inter-related motivations behind Definition \ref{def:matrix_disp_cnvx_entropy}. The first motivation comes from the notion of \emph{intrinsic dimensional functional inequalities}.
Consider a flow $(\prob_t)$ which is (approximately) trivial along certain directions in space, that is, its evolution  (approximately) takes place on a subspace of low dimension $\ll\dd$. In such settings, the explicit dependence on the ambient dimension $\dd$ in the notion of dimensional displacement convexity of the entropy, formulated as the concavity of $e^{-\frac{\Hent(t)}{\dd}}$, renders this notion oblivious to the intrinsic dimension of the flow $(\prob_t)$. Consequently, functional inequalities which are derived from dimensional displacement convexity are \emph{dimensional} functional inequalities, in the sense that the ambient dimension $\dd$ appears explicitly in the inequalities. This dimensional feature is undesirable in high-dimensional settings. On the other hand, in many practical settings, there is a lower-dimensional manifold inside the high-dimensional ambient space to which the objects of interest (approximately) belong (e.g., the \emph{manifold hypothesis}). In order to capture this phenomenon one needs \emph{intrinsic dimensional} functional inequalities where the ambient dimension is absent and which scale like the dimension of the object at hand. Indeed, it will be shown in this work that matrix displacement convexity allows to derive such {intrinsic dimensional functional inequalities, which improve on their classical dimensional counterparts by capturing more refined structures of the flow---see Section \ref{subsubsec:diff_inq}. This is because controlling the matrix entropy, rather than just its trace, facilitates the analysis of the flow $(\prob_t)$ along different directions in space.

The second motivation behind  Definition \ref{def:matrix_disp_cnvx_entropy} comes back to the discussion of curvature notions. For flat spaces, the $\text{CD}(0,\dd)$ curvature-dimension condition does not capture the full curvature structure of the space. Indeed, the $\text{CD}(0,\dd)$ condition implies a zero lower bound on the Ricci tensor, but in flat space one knows that the full Riemann tensor vanishes. More generally, there are important classes of manifolds where information beyond lower bounds on the Ricci tensor is given. One such prominent class in differential geometry is the class of Einstein manifolds with lower bounds on the \emph{sectional} curvature (which includes the sphere and hyperbolic space). What is the correct notion of displacement convexity that captures this type of curvature information? This question was taken up in \cite{KettererMondino2018}, but Definition \ref{def:matrix_disp_cnvx_entropy} seems to provide an alternative route as will be further explained in Section \ref{subsubsec:matrix_diff}. \\

To conclude this section the first result of this paper is stated informally:
\begin{theorem*}[Theorem \ref{thm:main_entropy}]
Suppose $(\prob_t,\phase_t)$ is a \nice flow satisfying \eqref{eq:flow_intro} and assume that $\noise$ is real, $\pot_t$ is convex, $\mf$ is concave, and $\f$ is non-decreasing. Then, $\HentM(t)$ is matrix displacement convex.
\end{theorem*}

\subsubsection{Matrix differential inequalities}
\label{subsubsec:matrix_diff}
One classical way to deduce convexity is via differential inequalities. The most basic example is expressing the displacement convexity of $\Hent(t)$ along some flow via the differential inequality
\begin{equation}
\label{eq:entcvxinq_intro}
\partial_t\speed(t)=\partial_{tt}^2\Hent(t)\ge 0. 
\end{equation}
The \emph{dimensional} displacement convexity of $\Hent(t)$ is equivalent to the differential inequality for the entropy production,
\begin{equation}
\label{eq:traceinq_intro}
\partial_t\speed(t)\ge \frac{1}{\dd}\speed^2(t). 
\end{equation}
In particular, comparing \eqref{eq:entcvxinq_intro} and \eqref{eq:traceinq_intro} shows that dimensional displacement convexity is stronger than displacement convexity. It will be shown in this work that the \emph{matrix} displacement convexity of $\HentM(t)$ follows from the \emph{matrix} differential inequality for the entropy production \emph{matrix},
\begin{equation}
\label{eq:matrixeinq_intro}
\partial_t\speedM(t)\succeq \speedM^2(t). 
\end{equation}
The inequality \eqref{eq:matrixeinq_intro} is stronger than \eqref{eq:traceinq_intro} by the Cauchy-Schwarz inequality. More importantly, the ambient dimension $\dd$ is absent from \eqref{eq:matrixeinq_intro}, and having an inequality for the full matrix (rather than just the trace as in \eqref{eq:traceinq_intro}) allows to control each direction of space separately. 

The proof of the matrix displacement convexity of $\HentM(t)$ will follow by establishing \eqref{eq:matrixeinq_intro}. In fact, more powerful differential matrix inequalities will be established, which in turn imply new  intrinsic dimensional functional inequalities. To state these differential inequalities  define the (positive semidefinite) \emph{Fisher information matrix} associated to a flow $(\prob_t)$ as
\[
\FishM(t):=\int_{\domain}(\nabla\log\prob_t)^{\otimes 2}\diff\prob_t
\]
and let 
\[
\costM_{\pm}(t):=\speedM(t)\pm\frac{\noise}{2}\FishM(t).
\]
\begin{theorem*}[Theorem \ref{thm:main_cost}, Theorem \ref{thm:main_entropy}]
Suppose $(\prob_t,\phase_t)$ is a \nice flow satisfying \eqref{eq:flow_intro} and assume  that $\noise$ is real, $\pot_t$ is convex, $\mf$ is concave, and $\f$ is non-decreasing.  Then,
\begin{equation}
\label{eq:cost_matrix_diff_inq_intro}
\partial_t\costM_{\pm}(t)\succeq \costM_{\pm}^2(t)+\int_{\domain} \nabla^2\pot_t\diff \prob_t +\int_{\domain} (-\nabla^2\mf)\mfstar \prob_t \diff \prob_t+\int_{\domain}\f'(\prob_t)(\nabla\prob_t)^{\otimes 2}\succeq \costM_{\pm}^2(t).
\end{equation}
Consequently,
\begin{equation}
\label{eq:entropy_matrix_diff_inq_intro}
\partial_t\speedM(t)\succeq \speedM^2(t)+\frac{\noise^2}{4}\FishM^2(t)+\int_{\domain} \nabla^2\pot_t\diff \prob_t +\int_{\domain} (-\nabla^2\mf)\mfstar \prob_t \diff \prob_t+\int_{\domain}\f'(\prob_t)(\nabla\prob_t)^{\otimes 2}\succeq \speedM^2(t).
\end{equation}
\end{theorem*}
The proofs of the matrix differential inequalities will rely on integration by parts, which provides a more analytic approach to the study of \eqref{eq:flow_intro} rather than relying on probabilistic representations. The proofs of the matrix differential inequalities  also shed light on the question posed at the end of Section \ref{subsubsec:matrixdspcv_intdimfuinq}. Crucial to the proofs of integration by parts is the exchange of derivatives, which is permitted  in the flat case treated in this work. But in the manifold setting, such an exchange of derivatives causes curvature terms to appear. Since the requisite differential inequalities are for \emph{matrices}, it is \emph{sectional} curvature terms which appear, rather than just Ricci terms. Hence, the verification of matrix differential inequalities, and hence matrix displacement convexity and intrinsic dimensional functional inequalities, is intimately tied to curvature information that goes beyond the classical curvature-dimension conditions. A concrete manifestation of this phenomenon can be found in \cite{eskenazis2023intrinsic} where Eskenazis and the author proved matrix differential inequalities (using different techniques) for heat flows over spaces of constant curvature. These matrix differential inequalities led to Hamilton-type inequalities (which operate on the matrix level and require assumptions on the full  Riemann tensor) which improve on Li-Yau inequalities (which operate on the trace level and only require information on the Ricci tensor).  The reader is referred to  \cite{eskenazis2023intrinsic} for further discussion.

\subsubsection{Intrinsic dimensional functional inequalities} 
\label{subsubsec:diff_inq}
Once the matrix differential inequalities \eqref{eq:cost_matrix_diff_inq_intro} and \eqref{eq:entropy_matrix_diff_inq_intro} are in place one can use known techniques to deduce functional inequalities. As explained above, the \emph{matrical} nature of the differential inequalities leads to the replacement of the (ambient) dimensional functional inequalities by more refined inequalities which capture the intrinsic dimension of the flow $(\prob_t)$, and thus improve on many classical results. These functional inequalities do not apply in the generality of the flows discussed above but apply in many cases of interest. Since the focus of this work is on matrix displacement convexity and the associated matrix differential inequalities, only some  intrinsic dimensional functional inequalities will be proven to show the power of the method. The reader is referred to the appropriate references for background on the significance of these functional inequalities.

Section \ref{sec:fun_inq} contains the  intrinsic dimensional functional inequalities which are summarized as follows:
\newline 

\begin{itemize}
\item \textbf{Theorem \ref{thm:entrpy_bounds_funinq_sec}.} Intrinsic dimensional lower and upper bounds on the growth  of the entropy $\Hent(t)$ along flows satisfying \eqref{eq:flow_intro}.
\newline

\item \textbf{Theorem \ref{thm:turnpike}.} Intrinsic dimensional  turnpike properties via dissipation of Fisher information along \emph{viscous} flows satisfying \eqref{eq:flow_intro}. 
\newline 

\item \textbf{Theorem \ref{thm:costinq}}. Intrinsic dimensional lower and upper bounds on certain costs associated to the flow \eqref{eq:flow_intro} when $\pot_t$ is independent of $t$, $\mf=0$, and $\noise\neq 0$.  These cost inequalities can also be seen as a generalization of the intrinsic dimensional local logarithmic Sobolev inequalities (Remark \ref{rem:LSI}).
\newline 

\item \textbf{Theorem \ref{thm:lontime_cost_energy}.} Intrinsic dimensional  long time asymptotics for cost and energy along entropic interpolation flows.
\newline 

\item \textbf{Theorem \ref{thm:EVI}.}  Intrinsic dimensional   evolution variational inequalities along entropic interpolation flows.
\newline 

\item \textbf{Theorem \ref{thm:Entcost_contract}.} Intrinsic dimensional contraction of entropic cost along entropic interpolation flows.
\end{itemize}

\begin{remark}
This work focuses on the development of matrix displacement convexity, and consequently intrinsic dimensional inequalities, for  certain functionals (e.g., entropy) along flows of the form \eqref{eq:flow_intro} in flat spaces. The natural next step is to characterize functionals which are matrix displacement convex (analogous to  \cite{McCann1997}), and to investigate the extension of the results of this paper to curved spaces.
\end{remark}

\subsection{Examples of density flows}
\label{subsec:examples}
To conclude the introduction this section demonstrates the scope of density flows of the  form \eqref{eq:flow_intro} via a number of important examples. The first step is to note that the equations in  \eqref{eq:flow_intro} have a variational characterization as the Euler-Lagrange equations of the functional
\begin{equation}
\label{eq:functional_intro}
(\prob,\vel): \partial_t\prob_t+\divnab\left(\prob_t\vel_t\right)=0\quad \mapsto \quad\int_0^{\ttime}\int_{\domain} \left[\lagrangian(t,x,\vel_t)+\frac{\noise^2}{8}|\nabla\log\prob_t|^2+\F(\prob_t)+\frac{1}{2}\mf\mfstar \prob_t\right]\diff \prob_t\diff t.
\end{equation}
Here, the Lagrangian  $\lagrangian$ is given by
\begin{equation}
\label{eq:Lagrangian_V_intro}
\lagrangian(t,x,\genvec):=\frac{|\genvec|^2}{2}-\pot_t(x),\quad t\in [0,\ttime],~ x\in\domain,~\genvec\in \R^{\dd},
\end{equation}
where $U_t:\domain \to \R$ is a potential term. The term $\mf\mfstar \prob_t$ stands for the convolution of the density $\prob_t$ with a symmetric interaction potential $\mf:\domain\to \R$, and $\F:\R_{\ge 0}\to \R$ is such that $\f(r)=\F(r)+r\F'(r)$ where $\f :\R_{\ge 0}\to \R$.
\begin{example}[Planning problems]
\label{ex:planning}
Consider the boundary conditions in \eqref{eq:flow_intro} specifying $\prob_0$ and $\prob_{\ttime}$. The \emph{planning problem} seeks to find  the optimal density flow going from $\prob_0$ to $\prob_{\ttime}$, subject to the minimization of the cost given by \eqref{eq:functional_intro}. The optimal flow $(\prob_t,\phase_t)$ is given by the equations \eqref{eq:flow_intro}.

\textbf{Optimal transport \cite[\S 5.4]{villani2021topics}.} Taking  $\pot_t=\f=\mf=0$ and $\noise =0$ leads to $(\prob_t,\phase_t)$ being the optimal flow minimizing \eqref{eq:functional_intro},
\begin{equation}
\label{eq:OT_functional_intro}
\int_0^{\ttime}\int_{\domain} \frac{|\vel_t|^2}{2}\diff \prob_t\diff t,
\end{equation}
among all flows satisfying the continuity equation $\partial_t\prob_t+\divnab(\prob_t\vel_t)=0$ with boundary conditions $\prob_0$ and $\prob_{\ttime}$. The function $\phase_t$ evolves according to 
\begin{equation}
\label{eq:optimal_transpot_intro}
\partial_t\phase_t+\frac{1}{2}|\nabla\phase_t|^2=0.
\end{equation}
In this setting the flow $(\prob_t)$ is a geodesic in the Wasserstein space between $\prob_0$ and $\prob_{\ttime}$, which is the fluid mechanics formulation by Benamou-Brenier \cite{BenamouBrenier2000} of the optimal transport problem between $\prob_0$ and $\prob_{\ttime}$.

\textbf{Heat flow.} Taking $\pot_t=\f=\mf=0$ and $\noise \to\infty$ leads to the flow $(\prob_t)$ corresponding to the heat equation
\begin{equation}
\label{eq:heat_intro}
\partial_t\prob_t-\frac{1}{2}\Delta\prob_t=0
\end{equation} 
(with boundary term $\prob_{\ttime}$ adjusted appropriately).

\textbf{Entropic interpolation \cite[\S 4.5]{chen2021stochastic}.} Taking  $\pot_t=\f=\mf=0$ leads to $(\prob_t,\phase_t)$ being the optimal flow minimizing \eqref{eq:functional_intro},
\begin{equation}
\label{eq:EOT_functional_intro}
\int_0^{\ttime}\int_{\domain} \frac{|\vel_t|^2}{2}\diff \prob_t\diff t+\frac{\noise^2}{8}\int_0^{\ttime}\int_{\domain}|\nabla\log\prob_t|^2\diff\prob_t\diff t,
\end{equation}
among all flows satisfying the continuity equation $\partial_t\prob_t+\divnab(\prob_t\vel_t)=0$ with boundary conditions $\prob_0$ and $\prob_{\ttime}$. The function $\phase_t$ evolves according to 
\begin{equation}
\label{eq:entropic_intro}
\partial_t\phase_t+\frac{1}{2}|\nabla\phase_t|^2+\frac{\noise^2}{2}\frac{\Delta \prob_t^{1/2}}{\prob_t^{1/2}}=0.
\end{equation} 
The flow $(\prob_t)$ is the entropic interpolation between $\prob_0$ and $\prob_{\ttime}$, which is the same flow of the \textbf{Schr\"odinger bridge problem} and is the dynamic formulation of  \textbf{entropic optimal transport}. The above fluid dynamics formulation (or stochastic control formulation) of entropic interpolation is due to Chen-Georgiou-Pavon \cite{ChenGeorgiouPavon2016} and Gentil-L\'{e}onard-Ripani \cite{GentilLeondardRipani2017}. The entropic interpolation flow encapsulates both the optimal transport flow and the heat flow as the limits $\noise\to 0$ and $\noise \to\infty$, respectively. 

\textbf{Regularization of planning problems \cite{GMST2019}.} Taking $\pot_t=\mf=0$ and $\noise =0$  leads to $(\prob_t,\phase_t)$ being the optimal flow minimizing \eqref{eq:functional_intro},
\begin{equation}
\label{eq:EOT_functional_intro}
\int_0^{\ttime}\int_{\domain} \frac{|\vel_t|^2}{2}\diff \prob_t\diff t+\int_0^{\ttime}\int_{\domain}\F(\prob_t)\diff\prob_t\diff t,
\end{equation}
among all flows satisfying the continuity equation $\partial_t\prob_t+\divnab(\prob_t\vel_t)=0$ with boundary conditions $\prob_0$ and $\prob_{\ttime}$. The function $\f$ is seen as a \emph{regularization} term of optimal transport and it is often assumed to be non-decreasing (an assumption under which the results of this work apply). The choice 
\[
\F(r)=\sgn\log r
\]
 leads to the \textbf{entropic regularization of optimal transport} as investigated by Porretta \cite{Porretta2023}.
\end{example}

\begin{example}[Mean-field games]
\label{ex:MFG}
The theory of mean-field games was developed by Huang-Malhame-Caines \cite{HMC2006} and Lasry-Lions \cite{lasry2007mean} to describe Nash equilibrium type concepts for games with large populations of agents. To describe this set up define the Hamiltonian $\hamiltonian$ to be the Legendre transform of the Lagrangian $\lagrangian$ (in the $\genvec$ variable), so that
\begin{equation}
\label{eq:Hamiltonian_intro}
\hamiltonian(t,x,\moment)=\frac{|\moment|^2}{2}+\pot_t(x),
\end{equation}
and let 
\begin{equation}
\label{eq:value'_intro}
\valuefn'(x,\prob):=\f(\prob(x))+(\mf\mfstar \prob)(x)
\end{equation}
for a density $\prob$ over $\domain$. Let $(\prob_t,\phase_t)$ be a flow satisfying \eqref{eq:flow_intro} and let $u_t:=\phase-\frac{\noise}{2}\log\prob_t$. Then, the system \eqref{eq:flow_intro} is the combination  of a Fokker-Planck equation and a Hamilton-Jacobi equation,
\begin{equation}
\label{eq:HJ_intro}
\begin{cases}
\partial_t\prob_t(x)-\frac{\noise}{2}\Delta\prob_t(x)+\divnab\left(\prob_t(x)\partial_{\moment}\hamiltonian(t,x,\nabla u_t)\right)=0,\\
\partial_t u_t(x)+\frac{\noise}{2}\Delta u_t(x)+\hamiltonian(t,x,\nabla u_t)=\valuefn'(x,\prob_t),
\end{cases}
\end{equation}
which describes the following stochastic optimal control problem. Consider an infinite population of agents where each agent evolves its state $x_t$ according to the stochastic differential equation
\begin{equation}
\label{eq:agent_intro}
dx_t=\vel_t(x_t)dt+\sqrt{\noise}dB_t, \quad x_0:=x\in \domain
\end{equation}
where $\vel_t$ is the control chosen by the agent and $(B_t)$ is a standard Brownian motion in $\R^{\dd}$. The agent's goal is to minimize 
\begin{equation}
\label{eq:EOT_functional_intro}
\mathbb E\left[\int_0^{\ttime}\lagrangian(t,x_t,\vel_t)\diff t+\int_0^{\ttime}\valuefn'(x_t,\prob_t)\diff t+u_{\ttime}(x_{\ttime})\right]
\end{equation}
where $\prob_t(x)$ is the density describing the fraction of agents at state $x$ at time $t$, and $u_{\ttime}$ stands for the cost at the final state $x_{\ttime}$.  The first term in \eqref{eq:EOT_functional_intro} stands for the energy spent by the control $\vel_t$, and the second term in \eqref{eq:EOT_functional_intro}  accounts for the effect of rest of the population of agents. For example, the common assumption in mean-field games (and in this work) that $\f$ is non-decreasing models the agent's aversion to overcrowding. At equilibrium, each agent chooses its control optimally and the resulting density is $\prob_t$ which satisfies the first equation in \eqref{eq:HJ_intro}. Letting $u_t(x)$ stand for the expected cost  that will be incurred by an agent playing optimally, starting at time $t$ at  state $x$,  one can show that $u_t$ solves the second equation in \eqref{eq:HJ_intro}. From a different perspective, $(\prob_t,\nabla u_t)$ can be derived as the optimal solution to the problem of minimizing 
\begin{equation}
\label{eq:MFG_functional_intro}
\int_0^{\ttime}\int_{\domain}\lagrangian(t,x,\vel_t)\diff \prob_t(x)\diff t+\int_0^{\ttime}\int_{\domain}\valuefn(x,\prob_t)\diff\prob_t\diff t+\int_{\domain}u_{\ttime}\diff\prob_{\ttime},
\end{equation}
among all flows $(\prob_t,\vel_t)$ satisfying the continuity equation $\partial_t\prob_t+\divnab(\prob_t\vel_t)=0$ with boundary conditions $\prob_0$ and $u_{\ttime}$, where
\begin{equation}
\label{eq:value_intro}
\valuefn(x,\prob):=\F(\prob(x))+\frac{1}{2}(\mf\mfstar \prob)(x),
\end{equation}
so that $\valuefn'$ is the functional derivative of $\valuefn$ (recall $\f(r)=\F(r)+r\F'(r)$). The functional \eqref{eq:MFG_functional_intro} is exactly \eqref{eq:functional_intro} and indeed  \eqref{eq:HJ_intro} is exactly \eqref{eq:flow_intro}. Equations \eqref{eq:HJ_intro} constitute a \emph{second-order} mean-field game system when  $\noise>0$ and a \emph{first-order} mean-field game system when  $\noise=0$. Note that in contrast to the planning problem of Example \ref{ex:planning}, where the boundary conditions were $(\prob_0,\prob_{\ttime})$, in the mean-field game setting the boundary conditions are $(\prob_0,u_{\ttime})$.
\end{example}

\begin{example}[Barotropic fluids] Let $e:\R_{\ge 0}\to \R$ be the internal energy of a fluid and let $\pressure: \R\to \R$ be the pressure function given by  $\pressure(r):=e'(r)r^2$. Taking $\pot_t=\mf=0, ~\noise =0$, and setting $e=-\F$, turns \eqref{eq:flow_intro} (after spatial differentiation of the second equation) into the system of equations 
\begin{equation}
\label{eq:Eulerintro}
\begin{cases}
\partial_t\prob_t+\divnab\left(\prob_t\nabla\phase_t\right)=0,\\
\partial_t\nabla\phase_t+\nabla_{\nabla\phase_t}\nabla\phase_t+\frac{p'(\prob_t)}{\prob_t}\nabla\prob_t=0.
\end{cases}
\end{equation}
The system \eqref{eq:Eulerintro} describes the \emph{compressible} Euler equations\footnote{The \emph{incompressible} Euler equations have the additional constraint $\Delta\phase_t=0$.} where the  pressure depends only on the fluid density $\prob_t$, which renders the fluid \emph{barotropic} \cite[\S 4.3]{KhesinMisiolekModin2021}. Normally, the pressure should be a non-decreasing function of the density, which translates to $\f'\le 0$ as $-\f'(r)=\frac{\pressure'(r)}{r}$. Most of the results of this work only apply to the setting $\f'\ge 0$ (which is the relevant setting for the mean-field games of Example \ref{ex:MFG}, and in principle can also be used in the planning problem of Example \ref{ex:planning}). However, there are in fact systems of fluid equations where $\f'\ge 0$, namely the zero-viscosity limit of the \textbf{slender jet equation} where $\pot_t(x)=gx$ (with $g>0$ standing for gravity) and $\f(r)=-\gamma r^{-\frac{1}{2}}$ (with $\gamma>0$ standing for the surface tension coefficient) \cite{constantin2020active, constantin2020hierarchy}. Note that in contrast to Example \ref{ex:planning} and Example \ref{ex:MFG}, the boundary conditions here are usually $(\prob_0,\nabla\phase_0)$.
\end{example}

\begin{example}[Semi-classical limits of non-linear Schr\"odinger equations] 
\label{ex:schrodinger}
Consider the equation 
\begin{equation}
\label{eq:schrodinger}
i\hbar\partial_t\Psi_t+\frac{\hbar^2}{2m}\Delta \Psi_t=\pot_t\Psi_t-(\mf\mfstar |\Psi_t|^2)\Psi_t-\f(|\Psi_t|^2)\Psi_t
\end{equation}
where $\Psi$ is a complex-valued wave function, $i$ is the imaginary unit, $m$ is the mass, and $\hbar$ is the reduced Planck constant. When $\f=\mf=0$, Equation \eqref{eq:schrodinger} is the standard linear Schr\"odinger equation with potential $\pot_t$ (often independent of $t$). The interaction potential $\mf(x)$ is often a power law, an inverse power law, or a logarithm in the norm $|x|$, and the non-linearity $\f(r)$ is often a polynomial or a logarithm in $r$. The connection between \eqref{eq:schrodinger} and \eqref{eq:flow_intro} is via the \textbf{Madelung transform} \cite{vonRenesse2012}:
Using the representation 
\[
\Psi_t(x):=\prob_t^{1/2}(x)e^{i\frac{m}{\hbar}\phase_t(x)},
\]
and assuming $|\Psi_t(x)|^2>0$ for every $x\in\domain$ and $t\in [0,\ttime]$, the flow $(\prob_t,\phase_t)$ satisfy
\begin{align}
\label{eq:entinterpol_schrodinger}
\begin{cases}
\partial_t\prob_t+\nabla{\cdot}(\prob_t\nabla\phase_t)=0,\\
\partial_t\phase_t+\frac{1}{2}|\nabla\phase_t|^2-\frac{\hbar^2}{2}\frac{\Delta \prob_t^{1/2}}{\prob_t^{1/2}}+\pot_t-\mf\mfstar \prob_t-\f(\prob_t)=0,
\end{cases}
\end{align}
where the units are chosen so that $m=1$. Equations \eqref{eq:entinterpol_schrodinger} are exactly the same as equations  \eqref{eq:flow_intro} with the choice $\noise =i\hbar$. The term  $\frac{\Delta \prob_t^{1/2}}{\prob_t^{1/2}}$ is known as the (non-local) \emph{quantum pressure} or \emph{Bohm potential}. Most of the results in this paper only apply to the case where $\noise$ is real so they cannot apply as is to \eqref{eq:entinterpol_schrodinger}. However, when taking $\hbar\to 0$, i.e., taking the \textbf{semi-classical limit} of \eqref{eq:schrodinger}, equations \eqref{eq:entinterpol_schrodinger} formally reduce to equations \eqref{eq:flow_intro} with $\noise =0$. In particular, the results of this work apply (at least formally) whenever $\pot_t$ is convex,  $\mf$ is concave, and $\f'\ge 0$. These assumptions cover a number of semi-classical limits of interest: 

\textbf{Semi-classical limit of the linear Schr\"odinger equation with convex potential.} Take $\f=\mf=0$ and $\pot_t$ to be convex. 

\textbf{Semi-classical limit of \emph{focusing} non-linear Schr\"odinger equations.} Take $\pot_t=\mf=0$ and $\f'\ge 0$. Prominent examples are 
\[
\f(r)=\sgn r,\quad \sgn > 0;
\]
 the \textbf{semi-classical limit of the focusing \emph{cubic} non-linear Schr\"odinger equation}, and 
 \[
 \f(r)=\sgn\log r,\quad \sgn > 0;
 \]
the \textbf{semi-classical limit of the focusing \emph{logarithmic} non-linear Schr\"odinger equation} (cf. entropic regularization of optimal transport in Example \ref{ex:planning}). 

For further information on semi-classical limits of non-linear focusing Schr\"odinger equations see \cite{McLaughlinMillerbook,miller1998focusingNLS,Grenier1998,ClarkeMiller02,Carles08} and \cite{bialynicki1976nonlinear,Pietro2014,Ferriere2020,carles2022logarithmic}.
\end{example}

The results of this work apply to all of the above examples and, in addition, to generalization afforded by considering general flows of the form \eqref{eq:flow_intro}.

\begin{remark}[Quantum drift-diffusion]
\label{rem:qunatum_diffusion_intro}
The \emph{quantum drift-diffusion model} \cite{quantumdrift-diffusion} (which for $\dd=1$ corresponds to the \emph{Derrida-Lebowitz-Speer-Spohn equation}) is defined by
\begin{align}
\label{eq:quantum drift-diffusion}
\begin{cases}
\partial_t\prob_t+\nabla{\cdot}(\prob_t\nabla\phase_t)=0,\\
\phase_t-2\frac{\Delta \prob_t^{1/2}}{\prob_t^{1/2}}=0.
\end{cases}
\end{align}
 It was shown by Gianazza-Savar\'e-Toscani  \cite{quantumdrift-diffusion} that the flow $(\prob_t,\phase_t)$ satisfying \eqref{eq:quantum drift-diffusion} is the gradient flow in Wasserstein space of the Fisher information functional\footnote{In Otto's seminal work \cite{Otto2001} it was shown that the heat equation is the gradient flow in Wasserstein space of the entropy functional.}, and an important part of the solution theory of \eqref{eq:quantum drift-diffusion} is the monotonicity of the entropy $\Hent(t)$ and Fisher information $\Tr[\FishM(t)]$. In Remark \ref{rem:qunatum_diffusion_sec} it is observed that, in addition to the known monotonicity of $\Tr[\FishM(t)]$, there is also  monotonicity (in the positive semidefinite sense) for the Fisher information \emph{matrix} $\FishM(t)$. This observation is in line with the theme of this work but the flow \eqref{eq:quantum drift-diffusion} does not fall under the framework of \eqref{eq:flow_intro},  and the monotonicity of the Fisher information matrix can anyway  be easily deduced from \cite{quantumdrift-diffusion}, so this observation will not be elaborated beyond Remark \ref{rem:qunatum_diffusion_sec}.
\end{remark}

\subsection*{Organization of paper} 
Section \ref{sec:preliminaries} establishes the assumptions, notation, and definitions used in this work. Section \ref{sec:formulas} contains the derivation of the formulas for the time derivatives of various quantities along the density flows. Section \ref{sec:matrix_diff_inq} contains the main results of this work where the differential matrix inequalities and matrix displacement convexity are derived. Finally,  Section \ref{sec:fun_inq} contains the intrinsic dimensional functional inequalities.

\subsection*{Acknowledgements} Many thanks to Giovanni Conforti for helpful remarks on the manuscript.  Thanks also to 
Alexandros Eskenazis, Matthew Rosenzweig, and Gigliola Staffilani for their comments and suggestions on this work. 

This material is based upon work supported by the National
Science Foundation under Award Numbers 2002022 and DMS-2331920.

\section{Preliminaries}
\label{sec:preliminaries}
This section collects the assumptions, notation, and definitions used in this work.

\subsection{Assumptions}
The existence  and regularity theory of density flows of the form \eqref{eq:flow_intro} is highly dependent on the precise form of the partial differential equations and the boundary conditions. Such questions are orthogonal to the topic of this work so will not be addressed here. Instead, sufficient regularity will be assumed to justify the computations. In certain settings, the results of this work are completely rigorous, provided sufficient regularity on the boundary conditions is assumed, while in other settings the computations are formal. In order to avoid distracting from the main point of this work this distinction will not be emphasized.

\begin{definition}[Nice flows]
\label{def:nice}
$~$

\begin{enumerate}
\item The domain $\domain$ is assumed to be a convex subset of $\R^{\dd}$ with smooth boundary (if the domain is bounded). 
\item The functions $\prob_t(x)$ and $\phase_t(x)$ are classical solutions of \eqref{eq:flow_intro}, differentiable in $t$, twice-differentiable in $x$, and are finite.
\item Integration by parts without boundary terms is justified. This entails either fast-enough decay of the flow and its derivative at infinity (when $\domain=\R^{\dd}$), or appropriate boundary conditions when $\domain$ is bounded. A good example to keep in mind is when $\domain$ is a flat torus in $\R^{\dd}$. 
\item The density flow $\prob_t$ has a smooth density with respect to the Lebesgue measure on $\R^{\dd}$, is assumed to be strictly positive, and integrates to 1, $\int_{\domain}\diff \prob_t(x)=1$ for all $t\in [0,\ttime]$.
\item The exchange of derivatives and integration is permitted.
\end{enumerate}
\end{definition}

\subsection{Notation} An absolutely continuous probability measure  $\gen$ will often be associated with its density with respect to the Lebesgue measure so that $\diff \nu=\nu \diff x$. To alleviate the notation the domain of spatial integrals will be omitted, $\int:=\int_{\domain}$, and the Lebesgue measure will be omitted as well, $\int:= \int \diff x$. Often, the $x$ argument of various functions will be omitted, e.g., $\int \gen(x)=\int \gen$, while the time dependence will be kept, e.g., $\int \vel_t(x)\diff x=\int \vel_t$. The metric on $\R^{\dd}$ is taken to be the standard Euclidean metric, denoted by $\langle \cdot,\cdot\rangle$, with the associated norm $|\cdot |$. The coordinates of a vector $\genvec\in \R^{\dd}$ are denoted by upper scripts, $\genvec = (\genvec ^1,\ldots,\genvec^{\dd})$.  The symmetric tensor product $\otimes_S$ is given by $\genvec \otimes_S \genvec':=\frac{1}{2}[\genvec\otimes \genvec' + \genvec'\otimes \genvec]$ for $\genvec,\genvec'\in \R^{\dd}$ where $\genvec\otimes \genvec'$ is the standard tensor product.

Matrix quantities will be denoted by calligraphic fonts, e.g., $\mathcal M$, and their traces (scalar) will be denoted by regular fonts, e.g., $M=\Tr[\mathcal M]$. The $(i,j)$th entry of $\mathcal M$ is denoted $\mathcal M_{ij}$.  The transport of a matrix $\mathcal M$ is given by $\mathcal M^{\mathsf{T}}$. The identity matrix on $\R^{\dd}$ is denoted by $\Id$. The symbols $\succeq$ and $\preceq$ will stand for the semi-definite order and will be applied only to symmetric matrices.

Time derivatives will be denoted as $\partial_t$ and spatial derivatives will be denoted as $\partial_i:=\partial_{x_i}$ and $\partial_{ij}^2:=\partial_{x_ix_j}^2$, etc. The spatial gradient and Hessian are denoted $\nabla,\nabla^2$, respectively, and $\divnab$ stands for the divergence of vector fields. Given a vector field $\vel$ over $\R^{\dd}$ denote by $\nabla \vel$ the matrix defined by $(\nabla \vel)_{ij}=\partial_i\vel^j$ with $\vel=(\vel^1,\ldots,\vel^{\dd})$, and write $\partial_k\vel:=(\partial_k\vel^1,\ldots,\partial_k\vel^{\dd})$ for $k\in \{1,\ldots, \dd\}$. The first and second derivative of a function $\eta$ over an interval are denoted by $\eta',\eta''$.

The summation $\sum_{k=1}^{\dd}$ will often be written as $\sum_k$.

\subsection{Definitions}
In the following definitions it is implicitly assumed that the expressions are well-defined. Throughout $\gen$ is a density (nonnegative function with unit integral) over $\domain$. 

The \emph{differential entropy} of $\gen$ is defined as
\begin{equation}
\label{eq:ent_def}
\Hent(\gen):=\int_{\domain} \log \nu \diff\nu.
\end{equation}
The \emph{Fisher information matrix} of $\gen$ is the symmetric matrix defined as
\begin{equation}
\label{eq:FishM_def}
\FishM(\gen):=-\int_{\domain} \nabla^2\log \gen\, \diff \gen = \int_{\domain} (\nabla\log \gen)^{\otimes 2}  \diff\nu,
\end{equation}
with the equality holding by integration by parts.
The \emph{Fisher information} of $\gen$ is 
\begin{equation}
\label{eq:Fish_def}
\Fish(\gen):=\Tr[\FishM(\gen)]=\int_{\domain}|\nabla\log \gen|^2\diff \gen.
\end{equation}
Let $(\prob_t)_{t\in [0,\ttime]}$ be a density flow which evolves according to a continuity equation:
\begin{equation}
\label{eq:continuity_entropysec}
\partial_t\prob_t+\divnab(\prob_t\vel_t)=0 \quad\forall~t\in [0,\ttime].
\end{equation}
For $t\in [0,\ttime]$ denote the entropy, Fisher information matrix, and Fisher information, respectively, of $\prob_t$ as
\begin{equation}
\label{eq:t_def}
\Hent(t):=\Hent(\prob_t),\quad \FishM(t):=\FishM(\prob_t),\quad \Fish(t):=\Tr[\FishM(t)].
\end{equation}
The \emph{entropy production matrix} is defined as
\begin{equation}
\label{eq:speedM_def}
\speedM(t):=\int_{\domain} (\nabla \prob_t\otimes_S  \vel_t)
\end{equation}
and the \emph{entropy production} is its trace
\begin{equation}
\label{eq:speed_def}
\speed(t):= \Tr[\speedM(t)]=\int_{\domain} \langle \nabla\prob_t ,\vel_t\rangle.
\end{equation}
The \emph{matrix entropy} is defined by
\begin{equation}
\label{eq:HentM_def}
\HentM(t) :=\int_0^t\speedM(s)\diff s
\end{equation}
so that
\[
\Hent(t)=\Hent(0)+\Tr[\HentM(t)].
\]
Another matrix which will play an important role comes from the driving vector field,
\begin{equation}
\label{eq:phaseM_def}
\phaseM(t):=\int_{\domain}  (\vel_t)^{\otimes 2}\diff\prob_t
\end{equation}
with its trace
\begin{equation}
\label{eq:phaseMnorm_def}
\phaseMnorm(t):=\Tr[\phaseM(t)].
\end{equation}
Note that $\phaseM(t)$ is symmetric. Finally, the following combinations of matrices will be crucial: Given $\noise \ge 0$ let
\begin{align}
\label{eq:CM_def}
&\costM_+(t):=\speedM(t)+\frac{\noise}{2}\FishM(t),\\
&\costM_-(t):=\speedM(t)-\frac{\noise}{2}\FishM(t).
\end{align}
The interpretation of $\costM_{\pm}(t)$ will become clearer in the subsequent sections.

\section{Density flows}
\label{sec:formulas}
This section derives the evolution equations of key quantities (entropy, entropy production matrix, etc.) along a density flow $(\prob_t)_{t\in [0,\ttime]}$ satisfying the continuity equation
\begin{equation}
\label{eq:continuity_eq}
\partial_t\prob_t+\divnab(\prob_t\vel_t)=0.
\end{equation}
Lemma \ref{lem:entropy_1der} and Lemma \ref{lem:Fish_1der} describe the first derivative along the flow of the entropy $\Hent(t)$ and Fisher information matrix $\FishM(t)$, respectively, for general flows satisfying \eqref{eq:continuity_eq}. While these results hold for general flows satisfying a continuity equation, the focus of this paper is on density flows of the form \eqref{eq:flow_intro}. Using the identity
\begin{equation}
\label{eq:quntum_pot_identity}
4\frac{\Delta \prob_t^{1/2}}{\prob_t^{1/2}}=|\nabla\log\prob_t|^2+2\Delta\log\prob_t,
\end{equation}
the flow \eqref{eq:flow_intro} can be written as
\begin{equation}
\label{eq:continuity_eq_phase}
\partial_t\prob_t+\divnab(\prob_t\nabla\phase_t)=0,
\end{equation}
with
\begin{equation}
\label{eq:phase_pde}
\partial_t\phase_t+\frac{1}{2}|\nabla\phase_t|^2+\frac{\noise^2}{8}\left[|\nabla\log\prob_t|^2+2\Delta\log\prob_t\right]+\pot_t-\mf\mfstar\prob_t-\f(\prob_t)=0,\quad \noise \in \R_{\ge 0}\cup i\R_{\ge 0}.
\end{equation}
For the class of flows $(\prob_t,\phase_t)$ satisfying \eqref{eq:continuity_eq_phase}-\eqref{eq:phase_pde}, Lemma \ref{lem:entropy_2der_matrix} provides a formula for the first derivative along the flow of the entropy production matrix $\speedM(t)$ and, consequently, deduces in  Corollary \ref{cor:entropy_2der} a formula for the second derivative of the entropy along the flow. In addition,  Lemma \ref{lem:phaseM_1stder} describes the evolution of $\phaseM(t)$ along the flow, and the evolution of its trace $\phaseMnorm(t)$ is given in Corollary \ref{cor:phaseMnorm_1stder}.\\

The first result describes the time evolution of the entropy of $\prob_t$. 
\begin{lemma}[1st derivative of entropy $\Hent(t)$]
\label{lem:entropy_1der}
Suppose $(\prob_t,\vel_t)_{t\in [0,\ttime]}$ is a \nice flow satisfying  \eqref{eq:continuity_eq}. Then,
\[
\partial_t \Hent(t) = S(t).
\]
\end{lemma}

\begin{proof}
From the continuity equation \eqref{eq:continuity_eq} and integration by parts,
\begin{align*}
\partial_t \Hent(t) &= \int (\partial_t \log \prob_t) \prob_t +  \int ( \log \prob_t) \partial_t\prob_t =- \int \nabla{\cdot}(\prob_t\vel_t) - \int (\log \prob_t)\nabla{\cdot}(\prob_t\vel_t)\\
&=0 + \int \langle \nabla\log \prob_t,\vel_t\rangle \prob_t = \int \langle \nabla\prob_t,\vel_t\rangle \\
&=S(t).
\end{align*}
\end{proof}

Next, the time evolution of the Fisher information matrix of $\prob_t$ is derived.
\begin{lemma}[1st derivative of Fisher information matrix $\FishM(t)$]
\label{lem:Fish_1der}
Suppose $(\prob_t,\vel_t)_{t\in [0,\ttime]}$ is a \nice flow satisfying \eqref{eq:continuity_eq}. Then,
\[
\partial_t\FishM(t)=\int_{\domain} [\nabla \vel_t\nabla^2\log\prob_t]\diff \prob_t+\int_{\domain}[\nabla \vel_t\nabla^2\log\prob_t]^{\mathsf{T}} \diff \prob_t.
\]
\end{lemma}

\begin{proof}
Recall that 
\[
\FishM_{ij}(t)=\int (\partial_i\log\prob_t\partial_j\log\prob_t)\prob_t=\int \partial_i\log\prob_t\partial_j\prob_t
\]
so that, by exchanging derivatives, 
\[
\partial_t\FishM_{ij}(t)=\int \partial_i\left(\frac{\partial_t\prob_t}{\prob_t}\right)\partial_j\prob_t+\int \partial_i\log\prob_t\partial_j\partial_t\prob_t.
\]
For a  vector field $\genvec=(\genvec^1,\ldots,\genvec^{\dd})$ and $k=1,\ldots,\dd$ let $\partial_k\genvec:=(\partial_k\genvec^1,\ldots,\partial_k\genvec^{\dd})$. Using the continuity equation \eqref{eq:continuity_eq} and exchanging derivatives gives
\begin{align*}
\partial_t\FishM_{ij}(t)&=-\int \partial_i\left(\frac{\nabla{\cdot}(\prob_t\vel_t)}{\prob_t}\right)\partial_j\prob_t-\int \partial_i\log\prob_t\,\partial_j\nabla{\cdot}(\prob_t\vel_t)\\
&=-\int \frac{\partial_i\nabla{\cdot}(\prob_t\vel_t)}{\prob_t}\partial_j\prob_t+\int \frac{\nabla{\cdot}(\prob_t\vel_t)}{\prob_t^2}\partial_i\prob_t\partial_j\prob_t-\int \partial_i\log\prob_t\,\partial_j\nabla{\cdot}(\prob_t\vel_t)\\
&=-\int [\nabla{\cdot} \partial_i(\prob_t\vel_t)]\partial_j\log\prob_t+\int \nabla{\cdot}(\prob_t\vel_t)\partial_i\log\prob_t\partial_j\log\prob_t-\int \partial_i\log\prob_t\,\nabla{\cdot} \partial_j(\prob_t\vel_t)\\
&=-\sum_k\int  \partial_{ki}^2(\prob_t\vel_t^k)\partial_j\log\prob_t+\sum_k\int \partial_k(\prob_t\vel_t^k)\partial_i\log\prob_t\partial_j\log\prob_t-\sum_k\int \partial_i\log\prob_t\, \partial_{kj}^2(\prob_t\vel_t^k).
\end{align*}
Hence, by integration by parts and exchanging derivatives,
\begin{align*}
&\partial_t\FishM_{ij}(t)=\sum_k \int  \partial_{i}(\prob_t\vel_t^k)\partial_{kj}^2\log\prob_t-\sum_k\int \prob_t\vel_t^k[\partial_{ik}^2\log\prob_t\partial_j\log\prob_t+\partial_i\log\prob_t\partial_{jk}^2\log\prob_t]+\sum_k\int \partial_{ik}^2\log \prob_t\partial_j(\prob_t\vel_t^k)\\
=&\sum_k \int  [\vel_t^k\partial_{i}\prob_t+\prob_t\partial_i\vel_t^k]\partial_{kj}^2\log\prob_t-\sum_k\int \vel_t^k[\partial_{ik}^2\log\prob_t\partial_j\prob_t+\partial_i\prob_t\partial_{jk}^2\log\prob_t]+\sum_k\int \partial_{ik}^2\log \prob_t[\vel_t^k\partial_j\prob_t+\prob_t\partial_j\vel_t^k]\\
=&\int [\nabla^2\log\prob_t\vel_t]_j\partial_i\prob_t+\int [\nabla\vel_t\nabla^2\log\prob_t]_{ij}\diff\prob_t-\int [\nabla^2\log\prob_t\vel_t]_{i}\partial_j\prob_t-\int [\nabla^2\log\prob_t\vel_t]_{j}\partial_i\prob_t\\
&+\int [\nabla^2\log\prob_t\vel_t]_i\partial_j\prob_t +\int  [\nabla\vel_t\nabla^2\log\prob_t]_{ji}\diff\prob_t\\
=&\int [\nabla\vel_t\nabla^2\log\prob_t]_{ij}\,\diff\prob_t+\int [\nabla\vel_t\nabla^2\log\prob_t]_{ji}\diff\prob_t.
\end{align*}
\end{proof}
The remainder of the section focuses on flows $(\prob_t,\phase_t)_{t\in [0,\ttime]}$ satisfying \eqref{eq:continuity_eq_phase}-\eqref{eq:phase_pde}. Note that under the evolution \eqref{eq:continuity_eq_phase}-\eqref{eq:phase_pde}, the entropy production matrix $\speedM(t)$ is symmetric since integration by parts gives
\begin{equation}
\label{eq:speedM_integration_by_parts}
\speedM(t)=\int_{\domain} (\nabla\prob_t\otimes_S\nabla\phase_t)\diff x=-\int_{\domain}  \phase_t\nabla^2\prob_t \diff x=-\int_{\domain} \nabla^2\phase_t\diff \prob_t.
\end{equation}

\begin{lemma}[1st derivative of entropy production matrix $\speedM(t)$]
\label{lem:entropy_2der_matrix}
Suppose $(\prob_t,\phase_t)_{t\in [0,\ttime]}$ is a \nice flow satisfying \eqref{eq:continuity_eq_phase}-\eqref{eq:phase_pde}. Then,
\[
\partial_t \speedM(t)=\int_{\domain}  (\nabla^2\phase_t)^2\diff \prob_t+\frac{\noise^2}{4}\int_{\domain}  (\nabla^2\log \prob_t)^2 \diff \prob_t+\int_{\domain}  \nabla^2\pot_t\diff\prob_t+\int_{\domain}  (-\nabla^2\mf)\mfstar\prob_t \diff\prob_t+\int_{\domain}  \f'(\prob_t)\frac{(\nabla\prob_t)^{\otimes 2}}{\prob_t} \diff \prob_t.
\]
\end{lemma}
\begin{proof}
By definition
\begin{align*}
\partial_t\speedM_{ij}(t)&=\frac{1}{2}\partial_t\int \partial_i\prob_t\partial_j\phase_t+\frac{1}{2}\partial_t\int \partial_j\prob_t\partial_i\phase_t\\
&=\frac{1}{2}\int \partial_t(\partial_i\prob_t)\partial_j\phase_t+\frac{1}{2}\int \partial_i\prob_t\partial_t(\partial_j\phase_t)+\frac{1}{2}\int \partial_t(\partial_j\prob_t)\partial_i\phase_t+\frac{1}{2}\int \partial_j\prob_t\partial_t(\partial_i\phase_t)\\
&=\frac{1}{2}\left[A_{ij}+B_{ij}+A_{ji}+B_{ji}\right]
\end{align*}
where
\begin{align*}
A_{ij}:=\int \partial_t(\partial_i\prob_t)\partial_j\phase_t,\quad \quad B_{ij}:=\int \partial_i\prob_t\partial_t(\partial_j\phase_t).
\end{align*}
To compute $A_{ij}$ note that by \eqref{eq:continuity_eq_phase},
\[
\partial_t\partial_i\prob_t=-\partial_i \nabla{{\cdot}} (\prob_t\nabla\phase_t)=-\sum_k\partial_i[\partial_k\phase_t\partial_k\prob_t+\prob_t\partial_{kk}^2\phase_t]=-\sum_k[\partial_{ik}^2\phase_t\partial_k\prob_t+\partial_k\phase_t\partial_{ik}^2\prob_t+\partial_i\prob_t\partial_{kk}^2\phase_t+\prob_t\partial_{kik}^3\phase_t].
\]
Hence, by integration by parts,
\begin{align*}
A_{ij}=&-\int \sum_k[\partial_{ik}^2\phase_t\partial_k\prob_t+\partial_k\phase_t\partial_{ik}^2\prob_t+\partial_i\prob_t\partial_{kk}^2\phase_t+\prob_t\partial_{kik}^3\phase_t]\partial_j\phase_t\\
=&-\int [\nabla^2\phase_t\nabla\prob_t]_i\partial_j\phase_t-\int [\nabla^2\prob_t\nabla\phase_t]_i\partial_j\phase_t+\sum_k\int \partial_k[\partial_j\phase_t\partial_i\prob_t]\partial_k\phase_t+\sum_k\int \partial_k[\prob_t\partial_j\phase_t]\partial_{ik}^2\phase_t\\
=&-\int [\nabla^2\phase_t\nabla\prob_t]_i\partial_j\phase_t-\int [\nabla^2\prob_t\nabla\phase_t]_i\partial_j\phase_t+\sum_k\int \partial_{jk}^2\phase_t\partial_i\prob_t\partial_k\phase_t+\sum_k\int \partial_j\phase_t\partial_{ik}^2\prob_t\partial_k\phase_t\\\
&+\sum_k\int \partial_k\prob_t\partial_j\phase_t\partial_{ik}^2\phase_t+\sum_k\int \prob_t\partial_{jk}^2\phase_t\partial_{ik}^2\phase_t\\
=&-\int [\nabla^2\phase_t\nabla\prob_t]_i\partial_j\phase_t-\int [\nabla^2\prob_t\nabla\phase_t]_i\partial_j\phase_t+\int [\nabla^2\phase_t\nabla\phase_t]_j\partial_i\prob_t+\int [\nabla^2\prob_t\nabla\phase_t]_i\partial_j\phase_t\\
&+\int [\nabla^2\phase_t\nabla\prob_t]_i\partial_j\phase_t+\int (\nabla^2\phase_t)^2_{ij}\diff \prob_t\\
=&\int [\nabla^2\phase_t\nabla\phase_t]_j\partial_i\prob_t+\int (\nabla^2\phase_t)^2_{ij}\diff \prob_t.
\end{align*}
To compute $B_{ij}$ note that by  \eqref{eq:phase_pde},
\begin{align*}
\partial_t\partial_j\phase_t&=-\partial_j\left\{\frac{1}{2}|\nabla\phase_t|^2+\frac{\noise^2}{8}\left[|\nabla\log\prob_t|^2+2\Delta\log\prob_t\right]+\pot_t-\mf\mfstar \prob_t-\f(\prob_t)\right\}\\
&=-\sum_k\partial_k\phase_t\partial_{jk}^2\phase_t-\frac{\noise^2}{4}\sum_k\partial_k\log\prob_t\partial_{jk}^2\log\prob_t-\frac{\noise^2}{4}\sum_k\partial_{kjk}^3\log\prob_t-\partial_j\pot_t+(\partial_j\mf)\mfstar \prob_t+\f'(\prob_t)\partial_j\prob_t.
\end{align*}
Hence, by integration by parts,
\begin{align*}
B_{ij}=&-\sum_k\int \partial_k\phase_t\partial_{jk}^2\phase_t\partial_i\prob_t-\frac{\noise^2}{4}\sum_k\int\partial_k\log\prob_t\partial_{jk}^2\log\prob_t \partial_i\prob_t-\frac{\noise^2}{4}\sum_k\int\partial_{kjk}^3\log\prob_t\partial_i\prob_t\\
&-\int\partial_j\pot_t\partial_i\prob_t+\int [(\partial_j\mf)\mfstar \prob_t]\partial_i\prob_t+\int\f'(\prob_t)\partial_j\prob_t\partial_i\prob_t\\
=&-\int [\nabla^2\phase_t\nabla\phase_t]_j\partial_i\prob_t+\frac{\noise^2}{4}\sum_k\int [\partial_{ik}^2\log\prob_t\partial_{kj}^2\log\prob_t+\partial_k\log\prob_t\partial_{ijk}^3\log\prob_t]\prob_t-\frac{\noise^2}{4}\sum_k\int\partial_{ijk}^3\log\prob_t\partial_k\prob_t\\
&+\int \partial_{ij}^2 \pot_t\diff \prob_t-\int (\partial_{ij}^2\mf)\mfstar \prob_t\diff\prob_t-\int\f'(\prob_t)\partial_{ij}^2\prob_t\diff\prob_t-\int\f''(\prob_t)\partial_i\prob_t\partial_j\prob_t\diff\prob_t\\
=&-\int [\nabla^2\phase_t\nabla\phase_t]_j\partial_i\prob_t+\frac{\noise^2}{4}\int (\nabla^2\log\prob_t)^2_{ij}\diff\prob_t+\frac{\noise^2}{4}\sum_k\int \partial_{ijk}^3\log\prob_t\partial_k\prob_t-\frac{\noise^2}{4}\sum_k\int\partial_{ijk}^3\log\prob_t\partial_k\prob_t\\
&+\int \partial_{ij}^2 \pot_t\diff \prob_t-\int (\partial_{ij}^2\mf)\mfstar \prob_t\diff\prob_t-\int\f'(\prob_t)\partial_{ij}^2\prob_t\diff\prob_t-\int\f''(\prob_t)\partial_i\prob_t\partial_j\prob_t\diff\prob_t\\
=&-\int [\nabla^2\phase_t\nabla\phase_t]_j\partial_i\prob_t+\frac{\noise^2}{4}\int (\nabla^2\log\prob_t)^2_{ij}\diff\prob_t+\int \partial_{ij}^2 \pot_t\diff \prob_t-\int (\partial_{ij}^2\mf)\mfstar \prob_t\diff\prob_t\\
&-\int\f'(\prob_t)\partial_{ij}^2\prob_t\diff\prob_t-\int\f''(\prob_t)\partial_i\prob_t\partial_j\prob_t\diff\prob_t.
\end{align*}
It follows that
\begin{align*}
A_{ij}+B_{ij}=&\int [\nabla^2\phase_t\nabla\phase_t]_j\partial_i\prob_t+\int (\nabla^2\phase_t)^2_{ij}\diff \prob_t-\int [\nabla^2\phase_t\nabla\phase_t]_j\partial_i\prob_t+\frac{\noise^2}{4}\int (\nabla^2\log\prob_t)^2_{ij}\diff\prob_t\\
&+\int \partial_{ij}^2 \pot_t\diff \prob_t-\int (\partial_{ij}^2\mf)\mfstar \prob_t\diff\prob_t-\int\f'(\prob_t)\partial_{ij}^2\prob_t\diff\prob_t-\int\f''(\prob_t)\partial_i\prob_t\partial_j\prob_t\diff\prob_t\\
=&\int (\nabla^2\phase_t)^2_{ij}\diff \prob_t+\frac{\noise^2}{4}\int (\nabla^2\log\prob_t)^2_{ij}\diff\prob_t+\int (\nabla^2\pot_t)_{ij}\diff \prob_t-\int (\nabla^2\mf)_{ij}\mfstar\prob_t \diff\prob_t\\
&-\int\f'(\prob_t)(\nabla^2\prob_t)_{ij}\diff\prob_t-\int\f''(\prob_t)(\nabla\prob_t)^{\otimes 2}_{ij}\diff\prob_t.
\end{align*}
Analogous argument applies to $A_{ji}+B_{ji}$. Finally, by integration by parts,
\begin{align*}
-\int\f'(\prob_t)(\nabla^2\prob_t)_{ij}\diff\prob_t=\int \f''(\prob_t)(\nabla\prob_t)^{\otimes 2}_{ij}\diff\prob_t+\int \f'(\prob_t)(\nabla\prob_t)^{\otimes 2}_{ij} \diff x
\end{align*}
so
\begin{align*}
-\int\f'(\prob_t)(\nabla^2\prob_t)_{ij}\diff\prob_t-\int\f''(\prob_t)(\nabla\prob_t)^{\otimes 2}_{ij}\diff\prob_t=\int \f'(\prob_t)(\nabla\prob_t)^{\otimes 2}_{ij} \diff x,
\end{align*}
which completes the proof.
\end{proof}
Combining Lemma \ref{lem:entropy_1der} and Lemma \ref{lem:entropy_2der_matrix} yields:
\begin{corollary}[2nd derivative of entropy $\Hent(t)$]
\label{cor:entropy_2der}
Suppose $(\prob_t,\phase_t)_{t\in [0,\ttime]}$ is a \nice flow satisfying \eqref{eq:continuity_eq_phase}-\eqref{eq:phase_pde}. Then,
\[
\partial_{tt}^2\Hent(t)=\int_{\domain}  \Tr[(\nabla^2\phase_t)^2]\diff \prob_t+\frac{\noise^2}{4}\int_{\domain}  \Tr[(\nabla^2\log \prob_t)^2] \diff \prob_t+\int_{\domain}  \Delta\pot_t\diff\prob_t+\int_{\domain}  (-\Delta \mf)\mfstar\prob_t\diff \prob_t+\int_{\domain}  f'(\prob_t)\frac{|\nabla\prob_t|^2}{\prob_t}\diff \prob_t.
\]
\end{corollary}

\begin{lemma}[1st derivative of  $\phaseM(t)$]
\label{lem:phaseM_1stder}
Suppose $(\prob_t,\phase_t)_{t\in [0,\ttime]}$ satisfy \eqref{eq:continuity_eq_phase}-\eqref{eq:phase_pde}. Then,
\[
\partial_t\phaseM_{ij}(t)=\frac{\noise^2}{4}\partial_t\FishM_{ij}(t)+\int_{\domain} \left\{\pot_t-\mf\mfstar\prob_t-\f(\prob_t)\right\}(\partial_i[ \prob_t\partial_j\phase_t]+\partial_j[ \prob_t\partial_i\phase_t]).
\]
\end{lemma}
\begin{proof}
Recall that
\[
\phaseM(t)=\int (\nabla\phase_t)^{\otimes 2}\diff \prob_t
\]
so that
\begin{align*}
\partial_t\phaseM_{ij}(t)&=\partial_t\int (\partial_i\phase_t\partial_j\phase_t)\prob_t=\int (\partial_i\partial_t\phase_t)\partial_j\phase_t\diff \prob_t+\int (\partial_j\partial_t\phase_t)\partial_i\phase_t\diff \prob_t+\int \partial_i\phase_t\partial_j\phase_t\partial_t\prob_t\\
&=A_{ij}+A_{ji}+B_{ij}
\end{align*}
where
\begin{align*}
&A_{ij}:=\int (\partial_i\partial_t\phase_t)\partial_j\phase_t\diff \prob_t,\quad\quad B_{ij}:=\int \partial_i\phase_t\partial_j\phase_t\partial_t\prob_t.
\end{align*}
To compute $A_{ij}$ note that by \eqref{eq:phase_pde},
\begin{align*}
A_{ij}=&-\int \partial_i\left\{\frac{1}{2}|\nabla\phase_t|^2+\frac{\noise^2}{8}\left[|\nabla\log\prob_t|^2+2\Delta\log\prob_t\right]+\pot_t-\mf\mfstar\prob_t-\f(\prob_t)\right\}\partial_j\phase_t\diff \prob_t\\
=&-\sum_k\int\left\{\partial_{ik}^2\phase_t\partial_k\phase_t+\frac{\noise^2}{4}\partial_{ik}^2\log\prob_t\partial_k\log\prob_t+\frac{\noise^2}{4}\partial_{kik}^3\log\prob_t\right\}\partial_j\phase_t\diff \prob_t\\
&-\int \partial_i\left\{\pot_t-\mf\mfstar\prob_t-\f(\prob_t)\right\}\partial_j\phase_t\diff \prob_t.
\end{align*}
Hence, by integration by parts,
\begin{align*}
A_{ij}=&-\int [\nabla^2\phase_t\nabla\phase_t]_i\partial_j\phase_t\diff\prob_t-\frac{\noise^2}{4}\int [\nabla^2\log\prob_t\nabla\log\prob_t]_i\partial_j\phase_t\diff\prob_t+\frac{\noise^2}{4}\sum_k\int (\partial_{ik}^2\log\prob_t)[\prob_t\partial_{jk}^2\phase_t+\partial_j\phase_t\partial_k\prob_t]\\
&+\int\left\{\pot_t-\mf\mfstar\prob_t-\f(\prob_t)\right\}\partial_i[ \prob_t\partial_j\phase_t]\\
=&-\int [\nabla^2\phase_t\nabla\phase_t]_i\partial_j\phase_t\diff\prob_t-\frac{\noise^2}{4}\int [\nabla^2\log\prob_t\nabla\log\prob_t]_i\partial_j\phase_t\diff\prob_t+\frac{\noise^2}{4}\int [\nabla^2\log\prob_t\nabla^2\phase_t]_{ij} \diff\prob_t\\
&+\frac{\noise^2}{4}\int [\nabla^2\log\prob_t\nabla\log\prob_t]_i\partial_j\phase_t\diff\prob_t+\int\left\{\pot_t-\mf\mfstar\prob_t-\f(\prob_t)\right\}\partial_i[ \prob_t\partial_j\phase_t]\\
=&-\int [\nabla^2\phase_t\nabla\phase_t]_i\partial_j\phase_t\diff\prob_t+\frac{\noise^2}{4}\int [\nabla^2\log\prob_t\nabla^2\phase_t]_{ij} \diff\prob_t+\int\left\{\pot_t-\mf\mfstar\prob_t-\f(\prob_t)\right\}\partial_i[ \prob_t\partial_j\phase_t].
\end{align*}
Analogously,
\begin{align*}
A_{ji}=-\int [\nabla^2\phase_t\nabla\phase_t]_j\partial_i\phase_t\diff\prob_t+\frac{\noise^2}{4}\int [\nabla^2\log\prob_t\nabla^2\phase_t]_{ji} \diff\prob_t+\int\left\{\pot_t-\mf\mfstar\prob_t-\f(\prob_t)\right\}\partial_j[ \prob_t\partial_i\phase_t].
\end{align*}
To compute $B_{ij}$ note that by \eqref{eq:continuity_eq_phase} and integration by parts,
\begin{align*}
B_{ij}=-\int \partial_i\phase_t \partial_j\phase_t \nabla{\cdot} (\prob_t\nabla\phase_t)=\sum_k\int \partial_k(\partial_i\phase_t \partial_j\phase_t)\partial_k\phase_t\diff\prob_t=\int [\nabla^2\phase_t\nabla\phase_t]_i\partial_j\phase_t\diff\prob_t+\int [\nabla^2\phase_t\nabla\phase_t]_j\partial_i\phase_t\diff\prob_t.
\end{align*}
It follows that
\begin{align*}
&A_{ij}+A_{ji}+B_{ij}\\
&=\frac{\noise^2}{4}\int [\nabla^2\phase_t\nabla^2\log\prob_t+\nabla^2\log\prob_t\nabla^2\phase_t]_{ij} \diff \prob_t+\int\left\{\pot_t-\mf\mfstar\prob_t-\f(\prob_t)\right\}(\partial_i[ \prob_t\partial_j\phase_t]+\partial_j[ \prob_t\partial_i\phase_t]).
\end{align*}
The proof is complete by Lemma \ref{lem:Fish_1der}.
\end{proof}
Combining Lemma \ref{lem:phaseM_1stder}  and  \eqref{eq:continuity_eq_phase} yields:
\begin{corollary}[1st derivative of $\phaseMnorm(t)$]
\label{cor:phaseMnorm_1stder}
Suppose $(\prob_t,\phase_t)_{t\in [0,\ttime]}$ is a \nice flow satisfying \eqref{eq:continuity_eq_phase}-\eqref{eq:phase_pde}. Then,
\[
\partial_t\phaseMnorm(t)=\frac{\noise^2}{4}\partial_t\Fish(t)-2\int_{\domain}\left\{\pot_t-\mf\mfstar\prob_t-\f(\prob_t)\right\}\partial_t\prob_t.
\]
\end{corollary}

\section{Matrix differential inequalities and matrix displacement convexity}
\label{sec:matrix_diff_inq}
In this section the main matrix differential inequalitiesof this work are derived. The main result is Theorem \ref{thm:matrix_diff_ineq} which provides matrix differential inequalities for $[0,\ttime]\ni t\mapsto \costM_{\pm}(t)$, for any flow satisfying \eqref{eq:continuity_eq_phase}-\eqref{eq:phase_pde}, \emph{provided that $\sigma\in \R_{\ge 0}$}. From Theorem \ref{thm:matrix_diff_ineq} it is possible to deduce a matrix differential inequality for $\speedM(t)$, which is the content of Theorem \ref{thm:speedM_diff_ineq}. In Section \ref{subsec:technical_matrix_diff}, a few technical results are collected which show how to obtain bounds on matrices and deduce matrix displacement convexity from matrix differential inequalities. Finally, Section \ref{subsec:main_results} apply Theorem \ref{thm:matrix_diff_ineq} and Theorem \ref{thm:speedM_diff_ineq}, together with the results of Section \ref{subsec:technical_matrix_diff}, to flows of the form \eqref{eq:continuity_eq_phase}-\eqref{eq:phase_pde} under convexity constraints.

The following theorem is the main result of this section and is based on the formulas of Section \ref{sec:formulas}.
\begin{theorem}[Matrix differential inequalities for $\costM_{\pm}(t)$]
\label{thm:matrix_diff_ineq}
Suppose $(\prob_t,\phase_t)_{t\in [0,\ttime]}$ is a \nice flow satisfying \eqref{eq:continuity_eq_phase}-\eqref{eq:phase_pde} with  $\noise \in \R_{\ge 0}$. Then,
\[
\partial_t\costM_+(t)\succeq \costM_+^2(t)+\int_{\domain} \nabla^2\pot_t\diff \prob_t+\int_{\domain} (-\nabla^2\mf)\mfstar\prob_t \diff\prob_t+\int_{\domain} \f'(\prob_t)\frac{(\nabla\prob_t)^{\otimes 2}}{\prob_t} \diff \prob_t,
\]
and
\[
\partial_t\costM_-(t)\succeq \costM_-^2(t)+\int_{\domain} \nabla^2\pot_t\diff \prob_t+\int_{\domain} (-\nabla^2\mf)\mfstar\prob_t \diff\prob_t+\int_{\domain} \f'(\prob_t)\frac{(\nabla\prob_t)^{\otimes 2}}{\prob_t} \diff \prob_t.
\]
\end{theorem}

\begin{proof}
Fix $\noise \in \R_{\ge 0}$. By Lemma \ref{lem:Fish_1der} and Lemma \ref{lem:entropy_2der_matrix},
\begin{align*}
\partial_t\costM_{\pm}(t)=&\int (\nabla^2\phase_t)^2\diff \prob_t+\frac{\noise^2}{4}\int (\nabla^2\log \prob_t)^2 \diff \prob_t+\int \nabla^2\pot_t\diff \prob_t+\int(-\nabla^2\mf)\mfstar\prob_t \diff\prob_t+\int \f'(\prob_t)\frac{(\nabla\prob_t)^{\otimes 2}}{\prob_t} \diff \prob_t\\
&\pm\frac{\noise}{2}\int [\nabla^2\phase_t\nabla^2\log\prob_t+\nabla^2\log\prob_t\nabla^2\phase_t] \diff \prob_t\\
=&\int \left[\nabla^2\phase_t\pm \frac{\noise}{2}\nabla^2\log\prob_t\right]^2\diff \prob_t+\int \nabla^2\pot_t\diff \prob_t+\int(-\nabla^2\mf)\mfstar\prob_t \diff\prob_t+\int \f'(\prob_t)\frac{(\nabla\prob_t)^{\otimes 2}}{\prob_t} \diff \prob_t\\
\succeq& \left[ \int \left(\nabla^2\phase_t\pm \frac{\noise}{2}\nabla^2\log\prob_t\right)\diff \prob_t\right]^2+\int \nabla^2\pot_t\diff \prob_t+\int(-\nabla^2\mf)\mfstar\prob_t \diff\prob_t+\int \f'(\prob_t)\frac{(\nabla\prob_t)^{\otimes 2}}{\prob_t} \diff \prob_t\\
=&\costM_{\pm}^2(t)+\int \nabla^2\pot_t\diff \prob_t+\int(-\nabla^2\mf)\mfstar\prob_t \diff\prob_t+\int \f'(\prob_t)\frac{(\nabla\prob_t)^{\otimes 2}}{\prob_t} \diff \prob_t,
\end{align*}
where the inequality holds by Jensen's inequality.
\end{proof}

By combining the differential inequalities of Theorem \ref{thm:matrix_diff_ineq} the following result is deduced.
\begin{theorem}[Matrix differential inequalities for entropy production matrix $\speedM(t)$]
\label{thm:speedM_diff_ineq}
Suppose $(\prob_t,\phase_t)_{t\in [0,\ttime]}$ is a \nice flow satisfying \eqref{eq:continuity_eq_phase}-\eqref{eq:phase_pde} with  $\noise \in \R_{\ge 0}$. Then,
\[
\partial_t\speedM(t)\succeq \speedM^2(t)+\frac{\noise^2}{4}\FishM^2(t)+\int_{\domain}\nabla^2\pot_t\diff \prob_t+\int_{\domain}(-\nabla^2\mf)\mfstar\prob_t \diff\prob_t+\int _{\domain}\f'(\prob_t)\frac{(\nabla\prob_t)^{\otimes 2}}{\prob_t} \diff \prob_t.
\]
\end{theorem}
\begin{proof}
Since
\[
\speedM(t)=\frac{\costM_{+}(t)}{2}+\frac{\costM_{-}(t)}{2}
\]
it follows from Theorem \ref{thm:matrix_diff_ineq} that
\begin{align*}
\partial_t\speedM(t)\succeq \frac{\costM_{+}^2(t)}{2}+\frac{\costM_{-}^2(t)}{2}+\int \nabla^2\pot_t\diff \prob_t+\int(-\nabla^2\mf)\mfstar\prob_t \diff\prob_t+\int \f'(\prob_t)\frac{(\nabla\prob_t)^{\otimes 2}}{\prob_t} \diff \prob_t.
\end{align*}
The result follows as
\begin{align*}
\frac{\costM_{+}^2(t)}{2}+\frac{\costM_{-}^2(t)}{2}&=\frac{1}{2}\left[\speedM^2(t)+2\speedM(t)\otimes_S\frac{\noise}{2}\FishM(t)+\frac{\noise^2}{4}\FishM^2(t)\right]+\frac{1}{2}\left[\speedM^2(t)-2\speedM(t)\otimes_S\frac{\noise}{2}\FishM(t)+\frac{\noise^2}{4}\FishM^2(t)\right]\\
&=\speedM^2(t)+\frac{\noise^2}{4}\FishM^2(t).
\end{align*}
\end{proof}
Both Theorem \ref{thm:matrix_diff_ineq} and Theorem \ref{thm:speedM_diff_ineq}  provide differential inequalities of the form
\begin{equation}
\label{eq:generic_matrix_diff_strong}
\partial_t M(t)\succeq M^2(t)+\text{remainder term}.
\end{equation}
In Section \ref{subsec:main_results} it will be shown that in many flows of interest the remainder term is nonnegative (in a semidefinite sense), which means that \eqref{eq:generic_matrix_diff_strong} implies  differential inequalities of the form
\begin{equation}
\label{eq:generic_matrix_diff}
\partial_t M(t)\succeq M^2(t).
\end{equation}
The following section shows how to take differential inequalities of the form \eqref{eq:generic_matrix_diff} and deduce bounds on $M(t)$ as well as obtain matrix displacement convexity.

\subsection{Matrix differential inequalities and displacement convexity} 
\label{subsec:technical_matrix_diff}
Suppose  for the rest of this section that $[0,\ttime]\ni t\mapsto M(t)$ is a differentiable function taking values in the set of $\dd\times \dd$ symmetric matrices. The first result shows how the differential inequality \eqref{eq:generic_matrix_diff} implies bounds on $M(t)$. 
\begin{lemma}
\label{lem:matrix_diff_to_scalar_diff}
 If 
\[
\partial_t M(t)\succeq M^2(t) \quad\forall ~t\in [0,\ttime],
\]
then, for any unit vector $\genvec\in \R^{\dd}$,
\begin{equation}
\label{eq:matrix_diff_to_scalar_diff}
\langle \genvec, M(t)\genvec\rangle \ge \frac{\langle \genvec, M(0)\genvec \rangle}{1-t\langle \genvec, M(0)\genvec \rangle} \quad\forall ~t\in [0,\ttime].
\end{equation}
\end{lemma}
\begin{proof}
Fix $\genvec\in \R^{\dd}$ and let $\eta(t):=\langle \genvec, M(t)\genvec\rangle$. Then, by the Cauchy-Schwarz inequality,
\begin{align*}
\partial_t\eta(t)=\langle \genvec, \partial_tM(t)\genvec\rangle\ge \langle \genvec, M^2(t)\genvec\rangle\ge \langle \genvec, M(t)\genvec\rangle^2=\eta^2(t). 
\end{align*}
The solution of the ordinary differential equation 
\[
\partial_t\xi(t)=\xi^2(t)\quad\forall~ t\in [0,\ttime]\quad \text{with}\quad \xi(0)=\eta(0)
\]
is $\xi(t):=\frac{\eta(0)}{1-t\eta(0)}$. Standard comparison \cite{Petrovitch} shows that $\eta(t)\ge \xi(t)$ for all $t\in [0,\ttime]$, which establishes \eqref{eq:matrix_diff_to_scalar_diff}. 
\end{proof}

\begin{corollary}
\label{cor:tr_inq}
Suppose
\[
\partial_t M(t)\succeq M^2(t) \quad\forall ~t\in [0,\ttime],
\]
and let $\{\lambda_i\}_{i=1}^{\dd}$ be the eigenvalues of $M(0)$. Then,
\[
\Tr[M(t)]\ge \sum_{i=1}^{\dd}\frac{\lambda_i}{1-\lambda_i t}.
\]
\end{corollary} 
\begin{proof}
Let $\{\genvec_i\}_{i=1}^{\dd}$ be the normalized eigenvectors of $M(0)$ corresponding to $\{\lambda_i\}_{i=1}^{\dd}$. By Lemma \ref{lem:matrix_diff_to_scalar_diff},
\[
\Tr[M(t)]=\sum_{i=1}^{\dd}\langle \genvec_i, M(t)\genvec_i\rangle \ge \sum_{i=1}^{\dd}\frac{\langle \genvec_i, M(0)\genvec_i \rangle}{1-t\langle \genvec_i, M(0)\genvec_i \rangle}=\sum_{i=1}^{\dd}\frac{\lambda_i}{1-\lambda_i t}.
\]
\end{proof}
Next it is shown how the differential inequality \eqref{eq:generic_matrix_diff} implies matrix  displacement convexity.
\begin{lemma}
\label{lem:refined_cnvx}
If
\[
\partial_t M(t)\succeq M^2(t) \quad\forall ~t\in [0,\ttime],
\]
then, $\int_0^tM(s)\diff s$ is matrix displacement convex, that is, for any unit vector $\genvec\in \R^{\dd}$, the function $\confn_{\genvec}:[0,\ttime]\to \R$ given by
\[
\confn_{\genvec}(t)=\exp\left[-\int_0^t\langle \genvec,M(s)\genvec\rangle\diff s\right]
\]
is concave. Consequently,
\begin{equation}
\label{eq:confn_prop1}
-\frac{1}{t}\le \langle \genvec,M(t)\genvec\rangle\le \frac{1}{\ttime-t}.
\end{equation}
\end{lemma}

\begin{proof}
To show the concavity of $\confn_{\genvec}$ it suffices to show that $\partial_{tt}^2\confn_{\genvec}(t)\le 0$ for every $t\in [0,\ttime]$. The first derivative is
\begin{align*}
\partial_t\confn_{\genvec}(t)=-\confn_{\genvec}(t)\langle \genvec,M(t)\genvec\rangle,
\end{align*}
and the second derivative is nonnegative as
\begin{align*}
\partial_{tt}^2\confn_{\genvec}(t)&=\confn_{\genvec}(t)\langle \genvec,M(t)\genvec\rangle^2-\confn_{\genvec}(t)\langle \genvec,\partial_tM(t)\genvec\rangle\\
&\le \confn_{\genvec}(t)\langle \genvec,M(t)\genvec\rangle^2-\confn_{\genvec}(t)\langle \genvec,M^2(t)\genvec\rangle\\
&=\confn_{\genvec}(t)\left\{\langle \genvec,M(t)\genvec\rangle^2-\langle \genvec,M^2(t)\genvec\rangle\right\}\\
&\le 0,
\end{align*}
where the first inequality holds by the assumption $\partial_t M(t)\succeq M^2(t)$, and the second inequality holds  by the Cauchy-Schwarz inequality. 

To establish \eqref{eq:confn_prop1} follow the argument of \cite[\S 3.3.1]{CCG22} and note that the concavity of $\confn_{\genvec}$ implies
\[
\partial_t\confn_{\genvec}(\ttime)\le \frac{\confn_{\genvec}(\ttime)-\confn_{\genvec}(t)}{\ttime -t}\le \partial_t\confn_{\genvec}(t)\le  \frac{\confn_{\genvec}(t)-\confn_{\genvec}(0)}{t}\le \partial_t\confn_{\genvec}(0).
\]
Since $\partial_t\confn_{\genvec}(t)=-\confn_{\genvec}(t)\langle \genvec,M(t)\genvec\rangle$,
and $\confn_{\genvec}(t)\ge 0$,
\[
-\confn_{\genvec}(t)\langle \genvec,M(t)\genvec\rangle=\partial_t\confn_{\genvec}(t)\le  \frac{\confn_{\genvec}(t)-\confn_{\genvec}(0)}{t}
\]
is equivalent to
\[
\langle \genvec,M(t)\genvec\rangle\ge  -\frac{1}{t}+\frac{\confn_{\genvec}(0)}{\confn_{\genvec}(t)t}\ge -\frac{1}{t}.
\]
The bound $\langle\genvec,M(t)\genvec\rangle\le \frac{1}{\ttime-t}$ follows analogously by using $\frac{\confn_{\genvec}(\ttime)-\confn_{\genvec}(t)}{\ttime -t}\le \partial_t\confn_{\genvec}(t)$.
\end{proof}

\begin{corollary}
\label{cor:tr_inq_intg}
Suppose
\[
\partial_t M(t)\succeq M^2(t) \quad\forall ~t\in [0,\ttime],
\]
and let $\{\lambda_i\}_{i=1}^{\dd}$ be the eigenvalues of $M(0)$. Then,
\[
\int_0^{\ttime}\Tr[M(t)]\diff t\ge -\sum_{i=1}^{\dd}\log\left(1-\ttime\lambda_i\right).
\]
\end{corollary} 
\begin{proof}
Taking $t=0$ in \eqref{eq:confn_prop1} gives $\lambda_i\le \frac{1}{\ttime}$. Hence,  $\lambda_i\le \frac{1}{t}$ for any $t\in [0,\ttime]$ which implies $0\le 1-t\lambda_i$ for any $t\in [0,\ttime]$. The result follows by integrating the bound in Corollary \ref{cor:tr_inq} from $t=0$ to $t=\ttime$.
\end{proof}

\subsection{Matrix differential inequalities and  displacement convexity along density flows}
\label{subsec:main_results}
This section shows that there are a number of important density flows where the matrix differential inequalities of Theorem \ref{thm:matrix_diff_ineq} and Theorem \ref{thm:speedM_diff_ineq} are of the form 
\begin{equation}
\label{eq:generic_matrix_diff_noneg_remainder}
\partial_t M(t)\succeq M^2(t)+\text{nonnegative term}.
\end{equation}
Hence, Lemma \ref{lem:matrix_diff_to_scalar_diff}, Corollary \ref{cor:tr_inq}, Lemma \ref{lem:refined_cnvx}, and Corollary \ref{cor:tr_inq_intg} are applicable.

\label{subsec:main_results}
\begin{theorem}[Differential inequalities  and matrix displacement convexity for $\costM_{\pm}(t)$]
\label{thm:main_cost}
Suppose $(\prob_t,\phase_t)_{t\in [0,\ttime]}$ is a \nice flow satisfying \eqref{eq:continuity_eq_phase}-\eqref{eq:phase_pde} with $\noise \in \R_{\ge 0}$, $\f'(r)\ge 0$ for every $r\in \R_{\ge 0}$, and $\int \{\nabla^2\pot_t-\nabla^2\mf\mfstar \prob_t\}\diff \prob_t \succeq 0$ for every $t\in [0,\ttime]$. Then, 
\begin{equation}
\label{eq:costM_diffinq_sec}
\partial_t\costM_{\pm}(t)\succeq \costM_{\pm}^2(t)+\int_{\domain} \nabla^2\pot_t\diff \prob_t +\int_{\domain} (-\nabla^2\mf)\mfstar \prob_t \diff \prob_t+\int_{\domain}\f'(\prob_t)\frac{(\nabla\prob_t)^{\otimes 2}}{\prob_t}\diff\prob_t\succeq \costM_{\pm}^2(t)\succeq 0.
\end{equation}
Consequently, for any unit vector $\genvec\in \R^{\dd}$,
\begin{equation}
\label{eq:thmTcostbound}
  \frac{\langle \genvec, \costM_{\pm}(0)\genvec \rangle}{1-t\langle \genvec, \costM_{\pm}(0)\genvec \rangle}\le \langle \genvec, \costM_{\pm}(t)\genvec\rangle \quad\forall ~t\in [0,\ttime],
\end{equation}
and
\begin{equation}
\label{eq:thmTcostboundTr}
 \sum_{i=1}^{\dd}\frac{\lambda_i}{1-\lambda_i t}\le \Tr[\costM_{\pm}(t)]\quad \text{where $\{\lambda_i\}_{i=1}^{\dd}$ are the eigenvalues of $\costM_{\pm}(0)$.}
\end{equation}
Furthermore, the matrix $\int_0^{t}\costM_{\pm}(s)\diff s$ is matrix displacement convex, that is, for any unit vector $\genvec\in \R^{\dd}$, the function $\confn_{\genvec}:[0,\ttime]\to \R$ given by
\begin{equation}
\label{eq:thmTrefined_cvx}
\confn_{\genvec}(t)=\exp\left[-\int_0^t\langle \genvec,\costM_{\pm}(s)\genvec\rangle\diff s\right] \text{ is concave}.
\end{equation}
Consequently, for every $t\in [0,\ttime]$,
\begin{equation}
\label{eq:thmTrefined_cvx_consq}
-\frac{1}{t}\le \langle \genvec, \costM_{\pm}(t)\genvec\rangle \le \frac{1}{\ttime -t},
\end{equation}
and
\begin{equation}
\label{eq:thmT_tr_inq_intg}
 -\sum_{i=1}^{\dd}\log\left(1-\ttime\lambda_i\right)\le \int_0^{\ttime}\Tr[\costM_{\pm}(t)]\diff t \quad \text{where $\{\lambda_i\}_{i=1}^{\dd}$ are the eigenvalues of $\costM_{\pm}(0)$.}
\end{equation}
\end{theorem}
The implications of Theorem \ref{thm:main_cost} to intrinsic dimensional functional inequalities will be derived in Section \ref{sec:fun_inq}.

The next result is analogous to Theorem \ref{thm:main_cost} but applies to $\speedM(t)$ rather than $\costM_{\pm}(t)$. Its implications to intrinsic dimensional functional inequalities will also be derived in Section \ref{sec:fun_inq}.
\begin{theorem}[Differential inequalities  and matrix displacement convexity for $\speedM(t)$]
\label{thm:main_entropy}
Suppose $(\prob_t,\phase_t)_{t\in [0,\ttime]}$ is a \nice flow satisfying \eqref{eq:continuity_eq_phase}-\eqref{eq:phase_pde} with $\noise \in \R_{\ge 0}$, $\f'(r)\ge 0$ for every $r\in \R_{\ge 0}$, and $\int \{\nabla^2\pot_t-\nabla^2\mf\mfstar \prob_t\}\diff \prob_t \succeq 0$ for every $t\in [0,\ttime]$. Then, 
\begin{equation}
\label{eq:S_diffinq_sec}
\partial_t\speedM(t)\succeq \speedM^2(t)+\frac{\noise^2}{4}\FishM^2(t)+\int_{\domain} \nabla^2\pot_t\diff \prob_t +\int_{\domain} (-\nabla^2\mf)\mfstar \prob_t \diff \prob_t+\int_{\domain}\f'(\prob_t)\frac{(\nabla\prob_t)^{\otimes 2}}{\prob_t}\diff \prob_t\succeq \speedM^2(t)\succeq 0.
\end{equation}
Consequently, for any unit vector $\genvec\in \R^{\dd}$,
\begin{equation}
\label{eq:thmEcostbound}
 \frac{\langle \genvec, \speedM(0)\genvec \rangle}{1-t\langle \genvec, \speedM(0)\genvec \rangle}\le \langle \genvec, \speedM(t)\genvec\rangle \quad\forall ~t\in [0,\ttime],
\end{equation}
and
\begin{equation}
\label{eq:thmEcostboundTr}
 \sum_{i=1}^{\dd}\frac{\lambda_i}{1-\lambda_i t}\le \speed(t)\quad\text{where $\{\lambda_i\}_{i=1}^{\dd}$ are the eigenvalues of $\speedM(0)$.}
\end{equation}
Furthermore, the matrix $\HentM(t)$ is matrix displacement convex, that is, for any unit vector $\genvec\in \R^{\dd}$, the function $\confn_{\genvec}:[0,\ttime]\to \R$ given by
\begin{equation}
\label{eq:thmErefined_cvx}
\confn_{\genvec}(t)=e^{-\langle \genvec,\HentM(t)\genvec\rangle} \text{ is concave}.
\end{equation}
Consequently, for every $t\in [0,\ttime]$
\begin{equation}
\label{eq:thmErefined_cvx_consq}
-\frac{1}{t}\le \langle \genvec, \speedM(t)\genvec\rangle \le \frac{1}{\ttime -t},
\end{equation}
and
\begin{equation}
\label{eq:thmE_tr_inq_intg}
-\sum_{i=1}^{\dd}\log\left(1-\ttime\lambda_i\right)\le\Tr[\HentM(\ttime)]=\Hent(\ttime)-\Hent(0) \quad \text{where $\{\lambda_i\}_{i=1}^{\dd}$ are the eigenvalues of $\speedM(0)$.}
\end{equation}
\end{theorem}

\begin{remark}[Quantum drift-diffusion]
\label{rem:qunatum_diffusion_sec} 
As mentioned in Remark \ref{rem:qunatum_diffusion_intro}, the quantum drift-diffusion model is given by the $(\prob_t,\phase_t)_{t\in [0,\ttime]}$ satisfying
\[
\begin{cases}
\partial_t\prob_t+\nabla{\cdot}(\prob_t\nabla\phase_t)=0,\\
\phase_t-2\frac{\Delta \prob_t^{1/2}}{\prob_t^{1/2}}=0.
\end{cases}
\]
In order to compute the derivatives of $\Hent(t)$ it turns out to be convenient to use the identity
\[
\nabla{\cdot}(\prob_t\nabla\phase_t)=\sum_{i,j}\partial_{ij}^2(\prob_t\partial_{ij}^2\log\prob_t).
\]
Then, the continuity equation implies (analogous to the proof of Lemma \ref{lem:entropy_1der}) that entropy decreases along the flow $(\prob_t)$ since
\[
\partial_t\Hent(t)=-\int\Tr[(\nabla^2\log\prob_t)^2]\diff \prob_t\le 0.
\]
For the computation of $\partial_t\FishM(t)$ apply Lemma \ref{lem:Fish_1der} to write
\[
\partial_t\FishM(t)=\int \left\{\nabla^2\left[2\frac{\Delta \prob_t^{1/2}}{\prob_t^{1/2}}\right]\nabla^2\log\prob_t+\nabla^2\log\prob_t\nabla^2\left[2\frac{\Delta \prob_t^{1/2}}{\prob_t^{1/2}}\right]\right\} \diff \prob_t,
\]
and use
\[
\nabla^2\log\prob_t=2\frac{\nabla^2\prob_t^{1/2}}{\prob_t^{1/2}}-2(\nabla\log \prob_t^{1/2})^{\otimes 2}
\]
to get
\begin{align*}
\int \nabla^2\left[2\frac{\Delta\prob_t^{1/2}}{\prob_t^{1/2}}\right]\nabla^2\log\prob_t\diff \prob_t&=\int \nabla^2\left[2\frac{\Delta\prob_t^{1/2}}{\prob_t^{1/2}}\right]2\frac{\nabla^2\prob_t^{1/2}}{\prob_t^{1/2}}\diff \prob_t-2\int \nabla^2\left[2\frac{\Delta\prob_t^{1/2}}{\prob_t^{1/2}}\right](\nabla\log\prob_t^{1/2})^{\otimes 2}\diff \prob_t\\
&=:A+B.
\end{align*}
Integration by parts shows that
\begin{align*}
A_{ij}=-\int\partial_{i}\left[2\frac{\Delta\prob_t^{1/2}}{\prob_t^{1/2}}\right]\partial_j\left[2\frac{\Delta\prob_t^{1/2}}{\prob_t^{1/2}}\right]\prob_t-2\int \partial_{i}\left[2\frac{\Delta\prob_t^{1/2}}{\prob_t^{1/2}}\right]\left\{(\nabla^2\prob_t^{1/2}\nabla\prob_t^{1/2})_j+\Delta\prob_t^{1/2}\partial_j\prob_t^{1/2}\right\},
\end{align*}
and
\begin{align*}
B_{ij}=2\int \partial_{i}\left[2\frac{\Delta\prob_t^{1/2}}{\prob_t^{1/2}}\right]\left\{\Delta\prob_t^{1/2}\partial_j\prob_t^{1/2}+(\nabla^2\prob_t^{1/2}\nabla \prob_t^{1/2})_j\right\}.
\end{align*}
It follows that
\[
A_{ij}+B_{ij}=-\int\partial_{i}\left[2\frac{\Delta\prob_t^{1/2}}{\prob_t^{1/2}}\right]\partial_j\left[2\frac{\Delta\prob_t^{1/2}}{\prob_t^{1/2}}\right]\prob_t
\]
and hence
\[
\partial_t\FishM(t)=-2\int \left(\nabla\left[2\frac{\Delta\prob_t^{1/2}}{\prob_t^{1/2}}\right]\right)^{\otimes 2}\diff\prob_t	\preceq 0,
\]
which establishes the monotonicity of the Fisher information \emph{matrix} along the quantum drift-diffusion flow.
\end{remark}

\section{Intrinsic dimensional functional inequalities}
\label{sec:fun_inq}
In this section Theorem \ref{thm:main_cost} and Theorem \ref{thm:main_entropy} will be used to derive intrinsic dimensional functional inequalities. When the boundary conditions of \eqref{eq:flow_intro} correspond to the planning problem, i.e., $(\prob_0,\prob_{\ttime})=(\sta,\fin)$ for densities $\sta,\fin$ over $\domain$, the time symmetry of the problem can be used:
\begin{remark}[Time symmetry]
\label{rem:time_symmetry}
The variational problem of \eqref{eq:functional_intro} with the boundary conditions  $(\sta,\fin)$ is time-symmetric. Consequently, if $(\prob_t,\phase_t)$ is the optimal flow with boundary conditions $(\sta,\fin)$, then the optimal flow with boundary conditions $(\fin,\sta)$ is  $(\tilde \prob_t,\tilde\phase_t)$ where $\tilde\prob_t:=\prob_{\ttime-t}$ and $\tilde \phase_t:=-\phase_{\ttime-t}$. Hence, the matrices $\tilde\costM_{\pm},\tilde\speedM,\tilde\FishM$ associated with  $(\tilde \prob_t,\tilde\phase_t)$  satisfy
\[
\tilde\speedM(t)=-\speedM(\ttime-t),\quad \tilde \FishM(t)=\FishM(\ttime-t),\quad \tilde\costM_{\pm}(t)=-\costM_{\mp}(\ttime -t)
\]
which implies
\[
\partial_t\tilde\costM_{\pm}(t)\succeq \tilde \costM_{\pm}^2(t),\quad \quad\partial_t\tilde\speedM(t)\succeq \tilde\speedM^2(t).
\]
\end{remark}

The first intrinsic dimensional functional inequality describes the growth of the entropy along the flow.

\begin{theorem}[Entropy growth]
\label{thm:entrpy_bounds_funinq_sec}
Suppose $(\prob_t,\phase_t)_{t\in [0,\ttime]}$ is a \nice flow satisfying \eqref{eq:continuity_eq_phase}-\eqref{eq:phase_pde} with $\noise \in \R_{\ge 0}$, $\f'(r)\ge 0$ for every $r\in \R_{\ge 0}$, and $\int \{\nabla^2\pot_t-\nabla^2\mf\mfstar \prob_t\}\diff \prob_t \succeq 0$ for every $t\in [0,\ttime]$. Then, 
\begin{equation}
\label{eq:entrpy_bounds_general}
 -\sum_{i=1}^{\dd}\log\left(1-\ttime\lambda_i(0)\right) \le \Hent(\ttime)-\Hent(0)
\end{equation}
where $\{\lambda_i(t)\}_{i=1}^{\dd}$ are the eigenvalues of $\speedM(t)$. Furthermore, under the planning problem boundary conditions $(\sta,\fin)$,
\begin{equation}
\label{eq:entrpy_bounds_planning}
 -\sum_{i=1}^{\dd}\log\left(1-\ttime\lambda_i(0)\right)\le \Hent(\ttime)-\Hent(0)\le  \sum_{i=1}^{\dd}\log\left(1+\ttime\lambda_i(\ttime)\right)
\end{equation}
where $\{\lambda_i(t)\}_{i=1}^{\dd}$ are the eigenvalues of $\speedM(t)$.
\end{theorem} 

\begin{proof}
Inequality \eqref{eq:entrpy_bounds_general} and the left-hand side of inequality \eqref{eq:entrpy_bounds_planning} is simply \eqref{eq:thmE_tr_inq_intg}. To get the right-hand side of inequality \eqref{eq:entrpy_bounds_planning}, Remark \ref{rem:time_symmetry} is used as follows. By \eqref{eq:thmEcostboundTr},
\begin{equation}
\label{eq:thm:entrpy_bounds_funinq_sec_proof}
 \sum_{i=1}^{\dd}\frac{\tilde\lambda_i(0)}{1-\tilde\lambda_i(0) t}\le \tilde\speed(t)
\end{equation}
where $\{\lambda_i(0)\}_{i=1}^{\dd}$ are the eigenvalues of $\tilde\speedM(0)=-\speedM(\ttime)$. By \eqref{eq:thmErefined_cvx_consq}, $\tilde\lambda_i(0)\le \frac{1}{\ttime}\le \frac{1}{t}$, which implies $1+t\lambda_i(\ttime)=1-t\tilde\lambda_i(0)\ge 0$. Hence, the integral over $t\in [0,\ttime]$ on the left-hand side of \eqref{eq:thm:entrpy_bounds_funinq_sec_proof} is equal to $-\sum_{i=1}^{\dd}\log\left(1+\ttime\lambda_i(\ttime)\right)$. The proof is complete by noting that $\int_0^{\ttime}\tilde\speed(t)\diff t=-\int_0^{\ttime}\speed(t)\diff t$.
\end{proof}

\subsection{Viscous flows} 
\label{subsec:viscous}
In this section the flow is assumed to be viscous, that is, $\noise \neq 0$. The first result pertains to the \emph{turnpike property}  of a  viscous flow $(\prob_t,\phase_t)_{t\in [0,\ttime]}$  satisfying \eqref{eq:continuity_eq_phase}-\eqref{eq:phase_pde}. The reader is referred to \cite{CCG22,faulwasser2022turnpike,turnpikeReview2022} for a discussion of the turnpike property, but in this context it suffices to state the formulation of the turnpike property by Clerc-Conforti-Gentil  \cite[Theorem 4.9]{CCG22}. They showed that when the flow $(\prob_t,\phase_t)$ is the \emph{entropic interpolation flow}, 
\begin{equation}
\label{eq:Thm4.9CCG22}
\Fish(t)\le \frac{n}{2t}+\frac{n}{2(\ttime -t)}.
\end{equation}
The next result improves on \eqref{eq:Thm4.9CCG22} by replacing the scalar inequality for the Fisher information by a matrix inequality for the Fisher information matrix, thus  disposing of the the ambient dimension $\dd$. In addition, the result applies to settings beyond entropic  interpolation.

\begin{theorem}[Turnpike properties via dissipation of Fisher information]
\label{thm:turnpike}
Suppose $(\prob_t,\phase_t)_{t\in [0,\ttime]}$ is a \nice flow satisfying \eqref{eq:continuity_eq_phase}-\eqref{eq:phase_pde} with $\noise >0$, $\f'(r)\ge 0$ for every $r\in \R_{\ge 0}$, and $\int \{\nabla^2\pot_t-\nabla^2\mf\mfstar \prob_t\}\diff \prob_t \succeq 0$ for every $t\in [0,\ttime]$. Then,
\[
 \FishM(t)\preceq \frac{1}{\noise}\left(\frac{1}{t}+\frac{1}{\ttime-t}\right)\Id.
\]
\end{theorem}

\begin{proof}
The proof is analogous to proof of   \cite[Theorem 4.9]{CCG22}. By \eqref{eq:thmTrefined_cvx_consq},
\[
\costM_+(t)\preceq \frac{1}{\ttime-t}\Id\quad\text{and}\quad  -\costM_-(t)\preceq \frac{1}{t}\Id
\]
so
\[
\noise\FishM(t)=\costM_+(t)-\costM_-(t)\preceq\left(\frac{1}{t}+\frac{1}{\ttime-t}\right)\Id.
\]
\end{proof}
The remainder of the results of this section are restricted to flows of the form 
\begin{equation}
\label{eq:viscous}
\begin{cases}
\partial_t\prob_t+\nabla{\cdot}(\prob_t\nabla\phase_t)=0, \quad \prob_0=\sta,~\prob_{\ttime}=\fin,\\
\partial_t\phase_t+\frac{1}{2}|\nabla\phase_t|^2+\frac{\noise^2}{8}\left[|\nabla\log\prob_t|^2+2\Delta\log\prob_t\right]+\pot-\f(\prob_t)=0, \quad\noise\neq 0,
\end{cases}
\end{equation}
so that the potential assumed to be independent of time, i.e.,
\[
\pot_t=\pot\quad\forall ~t\in [0,\ttime],
\]
and the interaction term $\mf$ is assumed to vanish. Under these assumptions an \emph{energy} can be defined which is constant along the flow. Begin by defining
\begin{align}
\label{eq:Hamilton_def}
\Hamilton(t)&:=\int_{\domain} \hamiltonian(x,\nabla\phase_t(x))\diff \prob_t(x) -\frac{\noise^2}{8}\Fish(t)-\int_{\domain} \F(\prob_t)\diff\prob_t\\
&=\int_{\domain}\left[\frac{1}{2}|\nabla\phase_t|^2+\pot-\frac{\noise^2}{8}|\nabla\log \prob_t|^2-\F(\prob_t)\right]\,\diff\prob_t,\nonumber
\end{align}
where the Hamiltonian $ \hamiltonian$ is given by
\[
 \hamiltonian(x,p):=\frac{|p|^2}{2}+\pot(x)
\]
and where $\F$ satisfies
\[
\f(r)=\F(r)+r\F'(r).
\]
\begin{lemma}[Preservation of energy]
\label{lem:Hamilton_preserve_trace}
Suppose $(\prob_t,\phase_t)_{t\in [0,\ttime]}$ is a \nice flow satisfying \eqref{eq:viscous}.  Then, the energy $\Hamilton(t)$ is constant along $t\in [0,\ttime]$. 
\end{lemma}

\begin{proof}
First note that
\[
\int_{\domain} \hamiltonian(x,\nabla\phase_t(x))\diff \prob_t(x) =\frac{1}{2}\phaseMnorm(t)+\int_{\domain} \pot \diff\prob_t.
\]
By Corollary \ref{cor:phaseMnorm_1stder},
\begin{align*}
\partial_t\frac{1}{2}\phaseMnorm(t)&=\frac{\noise^2}{8}\partial_t\Fish(t)-\int\left\{\pot-\f(\prob_t)\right\}\partial_t\prob_t
\end{align*}
while, on the other hand,
\begin{align*}
\partial_t\int \F(\prob_t)\diff\prob_t&=\int \prob_tF'(\prob_t)\partial_t\prob_t+\int F(\prob_t)\partial_t\prob_t=\int \f(\prob_t)\partial_t\prob_t.
\end{align*}
It follows that $\partial_t\Hamilton(t)=0$. 
\end{proof}
In light of Lemma \ref{lem:Hamilton_preserve_trace} define
\begin{equation}
\label{eq:Hamilton_trace_constant_def}
\Hamilton_{\ttime}:=\Hamilton(t)\quad\text{for any $t\in [0,\ttime]$}.
\end{equation}
Next define the \emph{cost}
\begin{align}
\label{eq:cost_fun_eq_sec}
\cost_{\ttime}&:=\int_0^{\ttime}\int_{\domain}\left[\lagrangian(x,\nabla\phase_t(x))+\frac{\noise^2}{8}|\nabla\log \prob_t(x)|^2+\F(\prob_t(x))\right]\,\diff\prob_t(x)\diff t\\
&=\int_0^{\ttime}\int_{\domain}\left[\frac{1}{2}|\nabla\phase_t|^2-\pot+\frac{\noise^2}{8}|\nabla\log \prob_t|^2+\F(\prob_t)\right]\,\diff\prob_t\diff t\nonumber,
\end{align}
where the Lagrangian $\lagrangian$ is given by
\[
\lagrangian(x,\genvec)=\frac{|\genvec|^2}{2}-\pot(x).
\]
The relation between the cost $\cost_{\ttime}$, the entropy $\Hent$, the energy $\Hamilton_{\ttime}$, and the matrix $\costM_{\pm}$  is captured by the following lemma.
\begin{lemma}
\label{lem:TrcostM}
Suppose $(\prob_t,\phase_t)_{t\in [0,\ttime]}$ is a \nice flow satisfying \eqref{eq:viscous}. Then,
\[
\frac{\noise}{2}\int_0^{\ttime}\Tr[\costM_{\pm}(t)]\diff t=\frac{\noise}{2}[\Hent(\ttime)-\Hent(0)]\pm[\cost_{\ttime}-\ttime \Hamilton_{\ttime}]\mp 2\int_0^{\ttime}\int_{\domain}[\F(\prob_t)-\pot]\diff\prob_t\diff t.
\]
\end{lemma}
\begin{proof}
By definition
\[
\cost_{\ttime}-\ttime \Hamilton_{\ttime}=\frac{\noise^2}{4}\int_0^{\ttime} \Fish(t)\diff t+2\int_0^{\ttime}\int[\F(\prob_t)-\pot]\diff\prob_t\diff t,
\]
so
\begin{align*}
\int_0^{\ttime} \Tr[\costM_{\pm}(t)]\diff t&=\int_0^{\ttime} \speed(t)\diff t\pm\frac{\noise}{2}\int_0^{\ttime} \Fish(t)\diff t=\Hent(\ttime)-\Hent(0)\pm\frac{\noise}{2}\int_0^{\ttime}\Fish(t)\diff t\\
&=\frac{2}{\noise}\left\{\frac{\noise}{2}[\Hent(\ttime)-\Hent(0)]\pm\frac{\noise^2}{4}\int_0^{\ttime}\Fish(t)\diff t\right\}\\
&=\frac{2}{\noise}\left\{\frac{\noise}{2}[\Hent(\ttime)-\Hent(0)]\pm[\cost_{\ttime}-\ttime \Hamilton_{\ttime}]\mp 2\int_0^{\ttime}\int[\F(\prob_t)-\pot]\diff\prob_t\diff t\right\}.
\end{align*}
\end{proof}
With Lemma \ref{lem:TrcostM} in hand the following intrinsic dimensional functional inequality for the combination of cost, entropy, and energy can be proved.

\begin{theorem}[Cost inequalities]
\label{thm:costinq}
Suppose $(\prob_t,\phase_t)_{t\in [0,\ttime]}$ is a \nice flow satisfying \eqref{eq:viscous} with $\noise >0$, $\f'(r)\ge 0$ for every $r\in \R_{\ge 0}$, and $\int \nabla^2\pot\diff \prob_t \succeq 0$ for every $t\in [0,\ttime]$. Then,
\begin{equation}
\label{eq:cost_bounds_general}
 -\frac{\noise}{2}\sum_{i=1}^{\dd}\log(1-\ttime\lambda_i(0))\le\frac{\noise}{2}[\Hent(\ttime)-\Hent(0)]\pm[\cost_{\ttime}-\ttime \Hamilton_{\ttime}]\mp 2\int_0^{\ttime}\int_{\domain}[\F(\prob_t)-\pot]\diff\prob_t\diff t
\end{equation}
where $\{\lambda_i(t)\}_{i=1}^{\dd}$ are the eigenvalues of $\costM_{\pm}(t)$. Furthermore, under the planning problem boundary conditions $(\sta,\fin)$,
\begin{align}
\label{eq:cost_bounds_planning}
 -\frac{\noise}{2}\sum_{i=1}^{\dd}\log(1-\ttime\lambda_i(0))\le\frac{\noise}{2}[\Hent(\ttime)-\Hent(0)]\pm[\cost_{\ttime}-\ttime \Hamilton_{\ttime}]\mp 2\int_0^{\ttime}\int_{\domain}[\F(\prob_t)-\pot]\diff\prob_t\diff t\le \frac{\noise}{2}\sum_{i=1}^{\dd}\log(1+\ttime\lambda_i(\ttime)),
\end{align}
where $\{\lambda_i(t)\}_{i=1}^{\dd}$ are the eigenvalues of $\costM_{\pm}(t)$.
\end{theorem}

\begin{proof}
Inequality \eqref{eq:cost_bounds_general} and the left-hand side of  inequality \eqref{eq:cost_bounds_planning} follows from \eqref{eq:thmT_tr_inq_intg} and Lemma \ref{lem:TrcostM}. For the right-hand side of inequality \eqref{eq:cost_bounds_planning}, use \eqref{eq:thmTcostboundTr} and Remark \ref{rem:time_symmetry} to get
\[
 \sum_{i=1}^{\dd}\frac{\lambda_i(0)}{1-t\lambda_i(0)}
\le \Tr[\costM_{\pm}(t)] \quad\text{and}\quad \sum_{i=1}^{\dd}\frac{\tilde\lambda_i(0)}{1-t\tilde\lambda_i(0)}\le \Tr[\tilde\costM_{\mp}(t)]=-\Tr[\costM_{\pm}(\ttime -t)],
\]
where $\{\lambda_i(0)\}_{i=1}^{\dd}$ are the eigenvalues of $\costM_{\pm}(0)$ and $\{\tilde\lambda_i(0)\}_{i=1}^{\dd}$ are the eigenvalues of $\tilde\costM_{\mp}(0)=-\costM_{\pm}(\ttime)$. Hence,
\[
 \sum_{i=1}^{\dd}\frac{\lambda_i(0)}{1-t\lambda_i(0)}
\le \Tr[\costM_{\pm}(t)]\quad\text{and}\quad \Tr[\costM_{\pm}(\ttime -t)]\le \sum_{i=1}^{\dd}\frac{\lambda_i(\ttime)}{1+t\lambda_i(\ttime)}.
\]
Using
\[
\int_0^{\ttime} \Tr[\costM_{\pm}(t)]\diff t=\int_0^{\ttime} \Tr[\costM_{\pm}(\ttime-t)]\diff t,
\]
and Lemma \ref{lem:TrcostM}, gives
\begin{equation}
\label{costinq_proof}
\frac{\noise}{2} \sum_{i=1}^{\dd}\int_0^{\ttime}\frac{\lambda_i(0)}{1-t\lambda_i(0)}\diff t\le\frac{\noise}{2}[\Hent(\fin)-\Hent(\sta)]\pm[\cost_{\ttime}-\ttime \Hamilton_{\ttime}]\mp 2\int_0^{\ttime}\int_{\domain}[\F(\prob_t)-\pot]\diff\prob_t\diff  t\le \frac{\noise}{2} \sum_{i=1}^{\dd}\int_0^{\ttime}\frac{\lambda_i(\ttime)}{1+t\lambda_i(\ttime)}\diff t.
\end{equation}
The next step is to note that by \eqref{eq:thmTrefined_cvx_consq}, $-\frac{1}{t}\le-\frac{1}{\ttime}\le  \lambda_i(\ttime)$, which implies that   $0\le  1+ t\lambda_i(\ttime)$. Hence, the integral over $t\in [0,\ttime]$ of the right-hand side of \eqref{costinq_proof} is equal to  $\frac{\noise}{2}\sum_{i=1}^{\dd}\log(1+\ttime\lambda_i(\ttime))$. Similarly, \eqref{eq:thmTrefined_cvx_consq} gives $\lambda_i(0)\le \frac{1}{\ttime}\le\frac{1}{t}$, which implies $1-t\lambda_i(0)\ge 0$.  Hence, the integral over $t\in [0,\ttime]$ of the left-hand side of \eqref{costinq_proof} is equal to $-\frac{\noise}{2}\sum_{i=1}^{\dd}\log(1-\ttime\lambda_i(0))$. 
\end{proof}

\begin{remark}[Intrinsic dimensional local logarithmic
Sobolev inequalities]
\label{rem:LSI}

The inequalities of Theorem \ref{thm:costinq} can be viewed as a generalization of the intrinsic dimensional  local  logarithmic
Sobolev inequalities for the Euclidean heat semigroup \cite[Equations (29) and (30)]{eskenazis2023intrinsic}. In particular, consider the entropic interpolation setting with  $\pot_t=\mf=\f=0$ and $\noise=1$. Fix $x\in \R^{\dd}$ and take $\sta:=\delta_x$ and $\diff\fin(y):=h(y)p_{\ttime}(x,y)$ where $h\ge 0$ and $p_{\ttime}$ is the heat kernel associated with the Euclidean  heat semigroup $\Pheat_{\ttime}$. Then, using the explicit expression for $(\prob_t)$ in \cite[Remark 4.1]{clerc2020variational}, one can formally derive the inequalities of \cite[Equations (29) and (30)]{eskenazis2023intrinsic} for the function $h$ evaluated at $x$:
\begin{equation}
\label{eq:intrinsicLSI}
\Pheat_{\ttime}(h\log h)-\Pheat_{\ttime}h\log \Pheat_{\ttime}h\le \frac{\ttime}{2} \Pheat_{\ttime}(\Delta h)+\frac{\Pheat_{\ttime}h}{2} \log\det\left\{\Id-\ttime \frac{\Pheat_{\ttime}(h\nabla^2\log h)}{\Pheat_{\ttime}h}\right\},
\end{equation}

\begin{equation}
\label{eq:reverseintrinsicLSI}
\frac{\ttime}{2} \Pheat_{\ttime}(\Delta h)-\frac{1}{2}\Pheat_{\ttime}h\log\det(\Id+\ttime \nabla^2\log \Pheat_{\ttime}h)\le \Pheat_{\ttime}(h\log h)-\Pheat_{\ttime}h\log \Pheat_{\ttime}h.
\end{equation}

 In \cite{clerc2020variational}, this entropic interpolation flow was used to prove the \emph{dimensional} local  log-Sobolev inequalities which are weaker than the \emph{intrinsic dimensional} local  log-Sobolev inequalities \eqref{eq:intrinsicLSI}-\eqref{eq:reverseintrinsicLSI}.
\end{remark}

\subsection{Entropic interpolation flows}
\label{subsec:fun_inq_particular}
In this section the flow will be assumed to be the entropic interpolation flow,
\begin{equation}
\label{eq:entropic_standard}
\begin{cases}
\partial_t\prob_t+\nabla{\cdot}(\prob_t\nabla\phase_t)=0, \quad \prob_0=\sta,~\prob_{\ttime}=\fin,\\
\partial_t\phase_t+\frac{1}{2}|\nabla\phase_t|^2+\frac{\noise^2}{8}\left[|\nabla\log\prob_t|^2+2\Delta\log\prob_t\right]=0, \quad\noise\in \R_{\ge 0},
\end{cases}
\end{equation}
over the domain $\domain=\R^{\dd}$. In contrast to Section \ref{subsec:viscous} where the energy preserved along the flow is a \emph{scalar} quantity, in the entropic interpolation setting a \emph{matrix} quantity is preserved as well. Define the \emph{matrix energy} as
\begin{equation}
\label{eq:HamiltonM_def}
\HamiltonM(t):=\frac{1}{2}\phaseM(t)-\frac{\noise^2}{8}\FishM(t).
\end{equation}

\begin{lemma}[Preservation of matrix energy]
\label{lem:preserve_HamiltonM}
Suppose that $(\prob_t,\phase_t)$ is a \nice flow satisfying \eqref{eq:entropic_standard}. Then, the energy matrix $\HamiltonM(t)$ is constant along $t\in [0,\ttime]$. 
\end{lemma}
\begin{proof}
The proof is immediate by Lemma \ref{lem:phaseM_1stder} which states
\[
\partial_t\phaseM(t)=\frac{\noise^2}{4}\partial_t\FishM(t).
\]
\end{proof}
In light of Lemma \ref{lem:preserve_HamiltonM} set
\begin{equation}
\label{eq:HamiltonM_con_def}
\HamiltonM_{\ttime}(\sta,\fin):=\HamiltonM(t)\quad\forall ~t\in [0,\ttime].
\end{equation}
In the setting of entropic interpolation a \emph{matrix} cost can be defined as
\begin{equation}
\label{eq:Cost_def_sec}
\Cost_{\ttime}(\sta,\fin)=\int_0^{\ttime}\int_{\R^{\dd}}\left[\frac{1}{2}(\nabla\phase_t)^{\otimes 2}+\frac{\noise^2}{8}(\nabla\log \prob_t)^{\otimes 2}\right]\,\diff\prob_t\diff t=\int_0^{\ttime}\left[\frac{1}{2}\phaseM(t)+\frac{\noise^2}{8}\FishM(t)\right]\diff t
\end{equation}
with its trace
\begin{equation}
\label{eq:Cost__scalar_def_sec}
\cost_{\ttime}(\sta,\fin):=\Tr[\Cost_{\ttime}(\sta,\fin)].
\end{equation}

The relation between the matrix cost and the matrix energy is captured by the following lemma, but first a remark is in order.
\begin{remark}[Time scaling]
\label{rem:time_scale}
Define
\begin{equation}
\label{eq:Acost_def}
\ACost_{\ttime}(\sta,\fin):=\inf_{(\prob_t,\vel_t)_{t\in [0,1]}}\int_0^1\int_{\R^{\dd}}\left[\frac{1}{2}(\vel_t)^{\otimes 2}+\ttime^2\frac{\noise^2}{8}(\nabla\log \prob_t)^{\otimes 2}\right]\,\diff\prob_t\diff t,
\end{equation}
where the minimum is over flows satisfying the continuity equations with boundary conditions:
\[
\partial_t\prob_t+\nabla{\cdot}(\prob_t\vel_t)=0 \quad\forall ~t\in [0,\ttime], \quad\prob_0=\sta,~\prob_1=\fin.
\]
Then, the Euler-Lagrange equations of \eqref{eq:Acost_def} are
\begin{equation}
\label{eq:entropic_time_rescale}
\begin{cases}
\partial_t\tilde\prob_t+\nabla{\cdot}(\tilde\prob_t\nabla\tilde\phase_t)=0, \quad \prob_0=\sta,~\prob_{\ttime}=\fin,\\
\partial_t\tilde\phase_t+\frac{1}{2}|\nabla\tilde\phase_t|^2+\ttime^2\frac{\noise^2}{8}\left[|\nabla\log\tilde\prob_t|^2+2\Delta\log\tilde\prob_t\right]=0,
\end{cases}
\end{equation}
and it is easy to see that $(\tilde\prob_t,\tilde\phase_t)=(\prob_{\ttime t},\ttime \phase_{\ttime t})$ where $(\prob_t,\phase_t)$ satisfy \eqref{eq:entropic_standard}.
\end{remark}
\begin{lemma}
\label{lem:cost_tau_der}
Suppose that $(\prob_t,\phase_t)$ is a \nice flow satisfying \eqref{eq:entropic_standard}. Assume that $\ttime\mapsto \partial_{\ttime}\Cost_{\ttime}(\sta,\fin)$ is differentiable. Then,
\[
\partial_{\ttime}\Cost_{\ttime}(\sta,\fin)=-\HamiltonM_{\ttime}(\sta,\fin).
\]
\end{lemma}

\begin{proof}
By the envelope theorem \cite{MilgromSegal2022} and Remark \ref{rem:time_scale},
\[
\partial_{\ttime}\ACost_{\ttime}(\sta,\fin)=\int_0^1\int_{\R^{\dd}}\ttime\frac{\noise^2}{4}(\nabla\log \tilde\prob_t)^{\otimes 2}\diff \tilde\prob_t
\]
where $(\tilde\prob_t)_{t\in [0,1]}$ is the optimal flow in $\ACost_{\ttime}(\sta,\fin)$ given by $\tilde\prob_t=\prob_{\ttime t}$. Hence, by the change of variables $t\mapsto \frac{t}{\ttime}$, 
\[
\partial_{\ttime}\ACost_{\ttime}(\sta,\fin)=\int_0^{\ttime}\int_{\R^{\dd}}\frac{\noise^2}{4}(\nabla\log \prob_t)^{\otimes 2}\diff \prob_t=\frac{\noise^2}{4}\int_0^{\ttime}\FishM(t)\diff t=\Cost_{\ttime}(\sta,\fin)-\ttime \HamiltonM_{\ttime}(\sta,\fin).
\]
On the other hand, changing variables $t\mapsto \ttime t$ shows that
\[
\ttime \Cost_{\ttime}(\sta,\fin)=\ACost_{\ttime}(\sta,\fin)
\] 
so it follows that
\[
\Cost_{\ttime}(\sta,\fin)-\ttime \HamiltonM_{\ttime}(\sta,\fin)=\partial_{\ttime}\ACost_{\ttime}(\sta,\fin)=\partial_{\ttime}[\ttime \Cost_{\ttime}(\sta,\fin)]=\Cost_{\ttime}(\sta,\fin)+\ttime\partial_{\ttime}\Cost_{\ttime}(\sta,\fin),
\]
which implies the result.
\end{proof}

To set up the first main result of this section recall the definition of the matrix entropy \eqref{eq:HentM_def} and
define
\begin{equation}
\label{eq:HentM_def_tau}
\HentM_{\ttime}(\sta,\fin):=\HentM(\ttime)
\end{equation}
so that
\begin{equation}
\label{eq:sta-fin-HentM_tau}
\Tr[\HentM_{\ttime}(\sta,\fin)]=\Hent(\fin)-\Hent(\sta).
\end{equation}

The following result is the intrinsic dimensional improvement of    \cite[Theorem 4.6]{CCG22} by Clerc-Conforti-Gentil, and is proved similarly.
\begin{theorem}[Large time asymptotics for cost and energy]
\label{thm:lontime_cost_energy}
Suppose that $(\prob_t,\phase_t)$ is a \nice flow satisfying \eqref{eq:entropic_standard}, and let $\ttime\ge 1$. Then,
\begin{equation}
\label{eq:HamiltonMw_bound}
-\HamiltonM_{\ttime}(\sta,\fin)\preceq \frac{\noise}{2}\frac{1}{\ttime}\Id,
\end{equation}
and, consequently,
\begin{equation}
\label{eq:CostMw_bound}
\Cost_{\ttime}(\sta,\fin)\preceq \Cost_1(\sta,\fin)+ \left(\frac{\noise}{2}\log\ttime\right) \Id.
\end{equation}
Moreover, for $t\in [0,\ttime)$,
\begin{equation}
\label{eq:phi_bound}
\int \left[\nabla\phase_t+\frac{\noise}{2}\nabla\log\prob_t\right]^{\otimes 2}\diff \prob_t \preceq \frac{\noise\HentM_{\ttime}(\sta,\fin)+2 \Cost_1(\sta,\fin)+ \left(\noise\log\ttime\right) \Id}{\ttime -t}.
\end{equation}

\end{theorem}

\begin{proof}
To prove \eqref{eq:HamiltonMw_bound} note that
\begin{equation}
\label{eq:phi}
0  \preceq \int \left[\nabla\phase_t+\frac{\noise}{2}\nabla\log\prob_t\right]^{\otimes 2}\diff \prob_t=\phaseM(t)+\noise \speedM(t)+\frac{\noise^2}{4}\FishM(t)
\end{equation}
which implies
\begin{equation}
\label{eq:HamiltonMw_bound_main_inq}
-\frac{1}{2}\phaseM(t)-\frac{1}{2}\frac{\noise^2}{4}\FishM(t)\preceq \frac{1}{2} \noise \speedM(t).
\end{equation}
Adding $\frac{\noise^2}{4}\FishM(t)$ to both sides of \eqref{eq:HamiltonMw_bound_main_inq} gives
\[
\frac{\noise^2}{8}\FishM(t)-\frac{1}{2}\phaseM(t)\preceq \frac{1}{2} \noise \speedM(t)+\frac{\noise^2}{4}\FishM(t)=\frac{\noise}{2}\costM_+(t)
\]
which is equivalent to
\[
-\HamiltonM_{\ttime}(\sta,\fin)\preceq\frac{\noise}{2} \costM_+(t).
\]
Taking $t=0$ gives
\[
-\HamiltonM_{\ttime}(\sta,\fin)\preceq\frac{\noise}{2}\costM_+(0) \preceq\frac{\noise}{2}\frac{1}{\ttime}\Id
\]
where the last inequality holds by \eqref{eq:thmTrefined_cvx_consq}. This establishes \eqref{eq:HamiltonMw_bound}.

To prove \eqref{eq:CostMw_bound} integrate \eqref{eq:HamiltonMw_bound} over $t$ from $0$ to $\ttime$ and use Lemma \ref{lem:cost_tau_der} to get
\begin{align*}
\Cost_{\ttime}(\sta,\fin)-\Cost_1(\sta,\fin)&=\int_1^{\ttime}\partial_s\Cost_{\ttime}(\sta,\fin) \diff s=\int_1^{\ttime}-\HamiltonM_{\ttime}(\sta,\fin) \diff s \\
&\preceq \frac{\noise}{2}\left(\int_1^{\ttime} \frac{1}{s}\diff s\right)\Id=\left(\frac{\noise}{2}\log\ttime\right)\Id.
\end{align*}
Finally, to prove \eqref{eq:phi_bound} note that by \eqref{eq:phi}, Lemma \ref{lem:phaseM_1stder}, and Theorem \ref{thm:matrix_diff_ineq},
\[
\partial_t \int \left[\nabla\phase_t+\frac{\noise}{2}\nabla\log\prob_t\right]^{\otimes 2}\diff \prob_t=\noise\partial_t\left[\speedM(t)+\frac{\noise}{2}\FishM(t)\right]=\noise \partial_t\costM_+(t) \succeq \noise\costM_+^2(t) \succeq 0,
\]
which shows that $[0,\ttime]\ni s\mapsto \int \left[\nabla\phase_s+\frac{\noise}{2}\nabla\log\prob_s\right]^{\otimes 2}\diff \prob_s$ is non-decreasing. Hence,
\begin{align*}
&(\ttime -t) \int \left[\nabla\phase_t+\frac{\noise}{2}\nabla\log\prob_t\right]^{\otimes 2}\diff \prob_t \preceq \int_t^{\ttime} \int \left[\nabla\phase_s+\frac{\noise}{2}\nabla\log\prob_s\right]^{\otimes 2}\diff \prob_s \diff s\\
&\preceq  \int_0^{\ttime} \int \left[\nabla\phase_s+\frac{\noise}{2}\nabla\log\prob_s\right]^{\otimes 2}\diff \prob_s \diff s=  \int_0^{\ttime} \left[\phaseM(s)+\noise \speedM(s)+\frac{\noise^2}{4}\FishM(s)\right] \diff s.
 \end{align*}
Since  
\[
 \int_0^{\ttime} \left[\phaseM(s)+\frac{\noise^2}{4}\FishM(s)\right]\diff s+\noise\int_0^{\ttime} \speedM(s)\diff s =2 \Cost_{\ttime}(\sta,\fin)+\noise\HentM_{\ttime}(\sta,\fin),
\]
the bound \eqref{eq:phi_bound}  follows by applying \eqref{eq:CostMw_bound}.
\end{proof}

The following result is the intrinsic dimensional improvement of the evolution variational inequality for the entropic cost of Ripani \cite[Corollary 11]{ripani2019convexity}, and is proved similarly.

\begin{theorem}[Evolution variational inequality]
\label{thm:EVI}
 Fix $\ttime =1, ~\noise=\sqrt{2}$, and suppose that $(\prob_t,\phase_t)$ is a \nice flow satisfying \eqref{eq:entropic_standard}. Let $(\Pheat_t)$ be the heat semigroup in $\R^{\dd}$ and suppose that  $t\mapsto \cost_1(\sta,\Pheat_t\fin)$ is differentiable. Then, for any normalized basis $\{\genvec_i\}_{i=1}^{\dd}$ of $\R^{\dd}$ and fixed $t\in [0,1]$,
\[
\partial_t\cost_{1}(\sta,\Pheat_t \fin)\le \frac{1}{2}\sum_{i=1}^{\dd}\left[1-e^{\langle \genvec_i,\HentM_{1}(\sta,\Pheat_t \fin)\genvec_i\rangle}\right].
\]
\end{theorem}

\begin{proof}
Let $\{\genvec_i\}_{i=1}^{\dd}$ be any normalized basis of $\R^{\dd}$ and fix $t\in [0,1]$. Consider the entropic interpolation flow betwen $\sta$ and $\Pheat_t\fin$ so, by \eqref{eq:thmErefined_cvx},  $\confn_{\genvec_i}(t):=e^{-\langle \genvec_i,\HentM(t)\genvec_i\rangle}$ is concave, and hence, 
\begin{equation}
\label{eq:conv_tangent}
\partial_t\confn_{\genvec_i}(1)\le \confn_{\genvec_i}(1)-\confn_{\genvec_i}(0).
\end{equation}
Inequality \eqref{eq:conv_tangent} is equivalent to
\[
-\confn_{\genvec_i}(1)\langle \genvec_i,\speedM(1)\genvec_i\rangle\le\confn_{\genvec_i}(1)-1,
\]
which upon rearrangement gives
\begin{equation}
\label{eq:conv_tangent_i}
-\langle\genvec_i,\speedM(1)\genvec_i\rangle\le 1-(\confn_{\genvec_i}(1))^{-1}.
\end{equation}
Using the definition of $\confn_{\genvec_i}(1)$, and summing over $i$ in \eqref{eq:conv_tangent_i}, gives 
\begin{align}
\label{eq:entropy_conv_tangent}
-\partial_t\Hent(t)|_{t=1}&=-\Tr[\speedM(1)]=-\sum_{i=1}^{\dd}\langle\genvec_i,\speedM(1)\genvec_i\rangle\le \sum_{i=1}^{\dd}\left[1-e^{\langle \genvec_i,\HentM(1)\genvec_i\rangle}\right]\nonumber\\
&=\sum_{i=1}^{\dd}\left[1-e^{\langle \genvec_i,\HentM_{1}(\sta,\Pheat_t \fin)\genvec_i\rangle}\right].
\end{align}
By \cite[Theorem 9]{ripani2019convexity},
\[
-\frac{1}{2}\partial_t\Hent(t)|_{t=1}=\partial_t\cost_{1}(\sta,\Pheat_t\fin)|_{t=0}
\]
so \eqref{eq:entropy_conv_tangent} implies
\begin{equation}
\label{eq:inq_at_0}
\partial_t\cost_{1}(\sta,\Pheat_t\fin)|_{t=0}\le \frac{1}{2}\sum_{i=1}^{\dd}\left[1-e^{\langle \genvec_i,\HentM_{1}(\sta,\Pheat_t \fin)\genvec_i\rangle}\right].
\end{equation}
By the semigroup property,  inequality \eqref{eq:inq_at_0} can be applied at any $t$ to yield the result.
\end{proof}

Finally, the last result is the  intrinsic dimensional improvement of  the entropic cost contraction of Ripani \cite[Corollary 13]{ripani2019convexity}, and is proved similarly.
\begin{theorem}[Contraction of entropy cost]
\label{thm:Entcost_contract}
 Fix $\ttime =1, ~\noise=\sqrt{2}$, and suppose that $(\prob_t,\phase_t)$ is a \nice flow satisfying \eqref{eq:entropic_standard}. Let $(\Pheat_t)$ be the heat semigroup in $\R^{\dd}$ and suppose that  $t\mapsto \cost_1(\sta,\Pheat_t\fin)$ is differentiable. Then,
\[
\cost_1(\Pheat_{\ttime}\sta,\Pheat_{\ttime} \fin)\le \cost_1(\sta, \fin)-2\sum_{i=1}^{\dd}\int_0^{\ttime}\sinh^2\left(\frac{\langle \genvec_i,\HentM_1(\Pheat_t \sta,\Pheat_t \fin)\genvec_i\rangle}{2}\right).
\]
\end{theorem}
\begin{proof}
Fix $s\in [0,1]$ and let $\{\genvec_i\}_{i=1}^{\dd}$ be any normalized basis of $\R^{\dd}$. Applying  Theorem \ref{thm:EVI} with $\sta\mapsto\Pheat_s\sta$ gives
\begin{equation}
\label{eq:s1}
\partial_t\cost_{1}(\Pheat_s\sta,\Pheat_t \fin)\le \frac{1}{2}\sum_{i=1}^{\dd}\left[1-e^{\langle \genvec_i,\HentM_{1}(\Pheat_s\sta,\Pheat_t \fin)\genvec_i\rangle}\right].
\end{equation}
On the other hand, by time symmetry (Remark \ref{rem:time_symmetry}),
\[
\partial_t\cost_{1}(\Pheat_t \fin,\sta)\le \frac{1}{2}\sum_{i=1}^{\dd}\left[1-e^{\langle \genvec_i,\HentM_{1}(\sta,\Pheat_t \fin)\genvec_i\rangle}\right],
\]
and switching the roles of $\sta$ and $\fin$ thus gives
\[
\partial_t\cost_{1}(\Pheat_t \sta,\fin)\le \frac{1}{2}\sum_{i=1}^{\dd}\left[1-e^{\langle \genvec_i,\HentM_{1}(\fin,\Pheat_t \sta)\genvec_i\rangle}\right].
\]
Taking $\fin\mapsto \Pheat_s \fin$ yields
\begin{equation}
\label{eq:s2}
\partial_t\cost_{1}(\Pheat_t \sta,\Pheat_s \fin)\le \frac{1}{2}\sum_{i=1}^{\dd}\left[1-e^{\langle \genvec_i,\HentM_{1}(\Pheat_s \fin,\Pheat_t \sta)\genvec_i\rangle}\right],
\end{equation}
and adding \eqref{eq:s1} and \eqref{eq:s2} shows that
\[
\partial_t\cost_{1}(\Pheat_s\sta,\Pheat_t \fin)+\partial_t\cost_{1}(\Pheat_t \sta,\Pheat_s \fin)\le \frac{1}{2}\sum_{i=1}^{\dd}\left[2-e^{\langle \genvec_i,\HentM_{1}(\Pheat_s\sta,\Pheat_t \fin)\genvec_i\rangle}-e^{\langle \genvec_i,\HentM_{1}(\Pheat_s \fin,\Pheat_t \sta)\genvec_i\rangle}\right].
\]
Taking $s=t$ yields
\begin{align*}
\partial_t\cost_{1}(\Pheat_t\sta,\Pheat_t \fin)&\le\frac{1}{2}\sum_{i=1}^{\dd}\left[2-e^{\langle \genvec_i,\HentM_{1}(\Pheat_t\sta,\Pheat_t \fin)\genvec_i\rangle}-e^{\langle \genvec_i,\HentM_{1}(\Pheat_t \fin,\Pheat_t \sta)\genvec_i\rangle}\right]\\
&=\sum_{i=1}^{\dd}\left[1-\left\{\frac{e^{\langle \genvec_i,\HentM_{1}(\Pheat_t\sta,\Pheat_t \fin)\genvec_i\rangle}+e^{-\langle \genvec_i,\HentM_{1}(\Pheat_t \sta,\Pheat_t \fin)\genvec_i\rangle}}{2}\right\}\right],
\end{align*}
where time  symmetry (Remark \ref{rem:time_symmetry}) was used to write $\tilde S(t)=-S(1-t)$ and hence $\HentM_{1}(\Pheat_t \fin,\Pheat_t \sta)=-\HentM_{1}(\Pheat_t \sta,\Pheat_t \fin)$. Integrating over $t$ from $0$ to $1$, and using $\cosh r=\frac{e^r+ e^{-r}}{2}$, gives
\begin{align*}
\cost_{1}(\Pheat_{\ttime}\sta,\Pheat_{\ttime} \fin)-\cost_{1}(\sta, \fin)=\sum_{i=1}^{\dd}\int_0^{\ttime}[1-\cosh(\langle \genvec_i,\HentM_{1}(\Pheat_t \sta,\Pheat_t \fin)\genvec_i\rangle)]\diff t.
\end{align*}
Finally, since $\sinh ^2 r=\frac{\cosh(2r)-1}{2}$, 
\begin{align*}
\cost_{1}(\Pheat_{\ttime}\sta,\Pheat_{\ttime} \fin)\le \cost_{1}(\sta, \fin)-2\sum_{i=1}^{\dd}\int_0^{\ttime}\sinh^2\left(\frac{\langle \genvec_i,\HentM_{1}(\Pheat_t \sta,\Pheat_t \fin)\genvec_i\rangle}{2}\right).
\end{align*}
\end{proof}

\bibliographystyle{amsplain0}
\bibliography{refEnt}

\end{document}